\documentclass[preprint,12pt]{elsarticle}

\usepackage[utf8]{inputenc}
\usepackage[T1]{fontenc}
\usepackage{amssymb,latexsym}
\usepackage{amsmath}
\usepackage{amsthm}
\usepackage{amssymb}
\usepackage{color}
\usepackage[shortlabels]{enumitem}
\usepackage[colorlinks=true,citecolor={blue},draft=false, urlcolor={black}]{hyperref}
\usepackage{cleveref}

\usepackage{ifdraft}
\usepackage{todonotes}

\usepackage{tikz}
\usetikzlibrary{calc,arrows,positioning,decorations.pathreplacing}

\def\Z{\mathbb Z}
\def\R{\mathbb R}
\def\Q{\mathbb Q}
\def\N{\mathbb N}
\def\A{\mathcal A}

\def\uu{\mathbf u}

\def\Z{\mathbb Z}
\def\R{\mathbb R}
\def\Q{\mathbb Q}
\def\N{\mathbb N}
\def\A{\mathcal A}

\def\Rcal{\mathcal R}
\def\M{\mathcal M}

\newcommand\oneint[8][]{
\pgfmathsetmacro\endX{(#2+#4)}
\begin{scope}[every path/.style={very thick},#1]
  \node(start#6)[at={(0,0)},label=left:{$0$}]{};
  \node(end#6)[at={(\endX,0)},label=right:{#8}]{};
  \draw[[-)] (start#6.center) -- (end#6.center);

  \ietpoint{discoAl}{#2}{#3}{alpha#6}

  \ietpoint{discoBe}{#4}{#5}{beta#6}

  #7
\end{scope}
}

\newcommand\ietpoint[5][pos=1,below]{
    \draw[style=#2] (#3,0.12) -- (#3,-0.12) node[#1]{#4};
    \node[ietinvisible] (#5) at (#3,0){};
}
\tikzset{ietinvisible/.style={outer sep=4,inner sep=0,minimum size=0}}
\newcounter{ietstep}
\newcounter{ietstepP}
\newcounter{ietstepPP}
\newcommand\ietstepcounters{   \stepcounter{ietstep}    \stepcounter{ietstepP}    \stepcounter{ietstepPP}}

\newtheorem{thm}{Theorem}
\newtheorem{theorem}[thm]{Theorem}
\newtheorem{coro}[thm]{Corollary}
\newtheorem{corollary}[thm]{Corollary}
\newtheorem{lem}[thm]{Lemma}

\newtheorem{claim}[thm]{Claim}

\newtheorem{proposition}[thm]{Proposition}
\newtheorem{prop}[thm]{Proposition}
\newtheorem{defi}[thm]{Definition}

\crefname{thm}{theorem}{theorems}
\crefname{theorem}{theorem}{theorems}
\crefname{coro}{corollary}{corollaries}
\crefname{example}{example}{examples}
\crefname{lem}{lemma}{lemmas}
\crefname{lmm}{lemma}{lemmas}
\crefname{claim}{claim}{claims}
\crefname{obs}{observation}{observations}
\crefname{proposition}{proposition}{propositions}
\crefname{prop}{proposition}{propositions}
\crefname{defi}{definition}{definitions}
\crefname{rem}{remark}{remarks}

\newtheorem{remark}[thm]{Remark}

\newtheorem{example}[thm]{Example}

\crefname{example}{example}{examples}

\begin{document}

\begin{frontmatter}

\title{On a faithful representation of Sturmian morphisms}

\author[bo,fn]{Jana Lep\v sov\'a}

\author[fn]{Edita  Pelantov\'a}

\author[fit]{\v St\v ep\'an Starosta\corref{ss}}

\ead{stepan.starosta@fit.cvut.cz}

\affiliation[bo]{organization={Univ. Bordeaux, CNRS, Bordeaux INP, LaBRI},
            addressline={UMR 5800},
            city={Bordeaux},
            postcode={F-33400},
            state={Talence},
            country={France}}

\affiliation[fn]{organization={Department of Mathematics, FNSPE, Czech Technical University in Prague},
            addressline={Trojanova 13},
            city={Prague},
            postcode={12000},
            country={Czech Republic}}

\affiliation[fit]{organization={Department of Applied Mathematics, FIT, Czech Technical University in Prague},
            addressline={Th\' akurova 9},
            city={Prague},
            postcode={16000},
            country={Czech Republic}}

\cortext[ss]{corresponding author}


\begin{abstract}
The set of morphisms mapping any Sturmian sequence to a Sturmian sequence forms together with composition the so-called monoid of Sturm.  For this monoid, we  define a   faithful representation  by  $(3\times 3)$-matrices with integer entries.
We find three convex cones in $\mathbb{R}^3$ and  show that a matrix  $R \in  Sl(\mathbb{Z},3)$ is a matrix  representing a Sturmian morphism if the  three cones are  invariant under multiplication by $R$ or $R^{-1}$.  This property offers a new tool to study Sturmian sequences.   We provide
alternative  proofs of four known results on Sturmian sequences fixed by a primitive morphism and a new result concerning the square root of a Sturmian sequence.
\end{abstract}


\begin{keyword}
Sturmian sequence \sep Sturmian morphism \sep square root of Sturmian sequence \sep faithful representation

\MSC 68R15
\end{keyword}

\end{frontmatter}




\section{Introduction}
Representation of objects  by linear operators is a  tool helping to  understand  the behaviour of algebraic structures, which   often serve as a model of  physical phenomena.  Representations of algebraic structures bring order into their description. Particle physics offers the most spectacular example:  the classification of elementary particles is based on an irreducible representation of certain Lie algebras \cite{algebra}.
In this article we demonstrate an application of the representation theory in combinatorics on words.
More specifically, we show that a representation of the special Sturmian  monoid provides a handy tool for proving some results on Sturmian sequences fixed by primitive morphisms.

  Sturmian sequences were introduced more than 80 years ago by Morse and Hedlund   \cite{HMo} and belong to the most explored objects in combinatorics on words.
In this article we focus on the Sturmian monoid, i.e. on the set of Sturmian morphisms, which are morphisms that  map any Sturmian sequence to  a Sturmian sequence. This set together with  composition of morphisms forms a monoid.
Association of each morphism with its incidence matrix forms a representation of this monoid. This representation maps the Sturmian monoid to a monoid of $ (2\times 2)$-matrices with integer entries. In~\Cref{FR},  we define  a new representation  $\mathcal{R}$ of the special Sturmian monoid  $\mathcal{M}$ by $ (3\times 3)$-matrices. Unlike  the representation by incidence matrices, our representation is faithful, i.e., $\mathcal{R}$  is injective. In~\Cref{submonoidSl}, we show that the matrices  assigned by $\mathcal{R}$ to Sturmian morphisms  form  a submonoid of  the group $Sl(\mathbb{Z},3)$ and  the  submonoid is characterized by  three convex cones in $\mathbb{R}^3$.
Using the  new representation $\mathcal{R}$, we provide new proofs of four known results in Section~\ref{app}.
A new result on square roots of Sturmian sequences, introduced in \cite{PeWh}, is obtained in~\Cref{ctverce}.
More specifically, we show that the fixing morphism of the square root of a Sturmian sequence $\uu$ can be found among the conjugates of small powers of the morphism fixing the $\uu$.

\section{Preliminaries}

Let $\mathcal{A}$ be an \emph{alphabet}, a finite set of \emph{letters}.
A \emph{(finite) word} $w$ is a finite sequence of elements of $\mathcal{A}$: $w = w_0w_1 \dots w_{n-1}$ with $w_i \in \mathcal{A}$.
The length of $w$, denoted $|w|$, equals $n$.
The \emph{empty word}, which is the unique word of length $0$, is denoted $\varepsilon$.
If $w = p f s$, i.e., the word $w$ is a concatenation of 3 words $p$, $f$ and $s$, we say that $p$ is a \emph{prefix} of $w$, $f$ is a \emph{factor} of $w$ and $s$ is a suffix of $w$.
The set of all words over $\mathcal{A}$ is denoted $\mathcal{A}^*$.

Let $\mathbf{s}= \left(s_i\right)_{i=0}^{+\infty}$ be a sequence over $\mathcal{A}$, that is, $s_i \in \mathcal{A}$ for all $i$.
Similarly, if $\mathbf{s} = pf
\mathbf{s'}$, where $p$ and $f$ are finite words and $\mathbf{s'}$ a sequence over $\mathcal{A}$, we say that $p$ is a \emph{prefix} of $\mathbf{s}$ and $f$ is a \emph{factor} of $\mathbf{s}$.
A factor $f$  of $\mathbf{s}$ is said to be \emph{left special} if  $af$ and $bf$ are also factors of $\mathbf{s}$ for  two distinct letters $a,b \in \mathcal{A}$.

The \emph{frequency of a letter $a \in \mathcal{A}$} in the sequence $\mathbf{s}$ equals $\lim_{i \to +\infty} \frac{| \textrm{pref}_i (\mathbf{s})|_a}{i}$, if the limit exists, where $\textrm{pref}_i (\mathbf{s})$ is the prefix of $\mathbf{s}$ of length $i$ and $|w|_a$ is the number of $a$'s in the word $w$.

A mapping $\mu: \mathcal{A}^* \to \mathcal{A}^*$ is a \emph{morphism} if for all $u,v \in \mathcal{A}^*$ we have $\mu(uv) = \mu(u)\mu(v)$.
A morphism $\mu$ is \emph{primitive} if there exists $k$ such that every $a \in \mathcal{A}$ is a factor of $\mu^k(b)$ for all $b \in \mathcal{A}$.

Morphisms naturally act on infinite sequences in the following manner.
Considering a sequence $\mathbf{s} = \left(s_i\right)_{i=0}^{+\infty}$ over $\mathcal{A}$, the image of $\mathbf{s}$ by the morphism $\mu$ is the sequence $\mu(\mathbf{s}) = \mu(s_0)\mu(s_1)\mu(s_2) \dots $ over $\mathcal{A}$.
We say that a sequence $\mathbf{s}$ is a \emph{fixed point} of $\mu$ if $\mu(\mathbf{s}) = \mathbf{s}$.

\section{Sturmian sequences}

Sturmian sequences allow many equivalent definitions, see for instance~\cite{Lo2} and~\cite{Fogg}.
We present a definition relying on mechanical sequences.
Let $\alpha \in [0,1]$ and $\delta$ be real numbers.
The sequences ${\bf s}_{\alpha, \delta}$  and ${\bf s}'_{\alpha, \delta}$ defined by
$${\bf s}_{\alpha, \delta}(n) := \lfloor \alpha(n+1) +\delta\rfloor -\lfloor \alpha n +\delta\rfloor   \quad \text{ for each } n \in \N $$
and
$${\bf s}'_{\alpha, \delta}(n) := \lceil\alpha(n+1) +\delta\rceil-\lceil \alpha n +\delta\rceil  \quad \text{ for each } n \in \N $$
are called the \emph{lower} and the \emph{upper}, respectively, \emph{mechanical sequences} with the \emph{slope} $\alpha$ and the \emph{intercept}~$\delta$.

It is easy to see that the sequences ${\bf s}_{\alpha, \delta}$  and ${\bf s}'_{\alpha, \delta}$ have elements in the set $\{0,1\}$.
If $\alpha$ is an irrational number, then ${\bf s}_{\alpha, \delta}(n)$ is a \emph{lower Sturmian sequence} and ${\bf s}'_{\alpha, \delta}(n)$
is an \emph{upper Sturmian sequence}.

Moreover, ${\bf s}_{\alpha, \delta}(n)=0$ if the fractional part of $\alpha n +\delta$ belongs to $ [0, 1-\alpha)$,  otherwise ${\bf s}_{\alpha, \delta}(n)=1$.
Similarly, ${\bf s'}_{\alpha, \delta}(n)=0$ if the fractional part of $\alpha n +\delta$ belongs to $ (0, 1-\alpha]$,  otherwise ${\bf s}_{\alpha, \delta}(n)=1$.
This property leads to an equivalent definition of Sturmian sequences, which relies on the two interval exchange transformation.

For given parameters $\ell_0, \ell_1 >0$, we consider two intervals of length $\ell_0$ and $\ell_1$.
To define a lower Sturmian sequence, we use the left-closed right-open intervals   $I_0=[0,\ell_0)$ and $I_1=[\ell_0, \ell_0+\ell_1)$, to define an upper Sturmian sequence, we use the left-open  right-closed intervals   $I_0=(0,\ell_0]$ and $I_1=(\ell_0, \ell_0+\ell_1]$.
The \emph{two interval exchange transformation} (2iet) $T: I_0\cup I_1\to I_0\cup I_1$ is defined by
$$
T(x) =
\left\{
\begin{array}{ll}
x + \ell_1 & \quad \text{if} \ x \in I_0,\\
x -  \ell_0 & \quad \text{if}  \ x \in I_1.
\end{array}
\right.
$$
If we take an initial point $ \rho \in I_0\cup I_1$, the sequence $\uu = u_0u_1u_2 \dots \in \{0,1\}^\N$ defined by
$$
u_n =
\left\{
\begin{array}{ll}
 0 & \quad \text{if} \ T^n(\rho) \in  I_0,\\
 1 & \quad \text{if} \ T^n(\rho)  \in I_1,
\end{array}
\right.
$$
i.e., a coding of the trajectory of the point $\rho$, is a \textit{2iet sequence} with the \textit{parameters} $\ell_0,\ell_1,\rho$.
We shall use the following notation for this fact:
$$ \vec{v}(\uu)= (\ell_0, \ell_1, \rho)^\top,$$
and refer to $(\ell_0, \ell_1, \rho)^\top$ as a \emph{vector of parameters of $\uu$}.

Clearly, collinear triples $(\ell_0, \ell_1, \rho)^\top$ and $c(\ell_0, \ell_1, \rho)^\top$ produce the same infinite sequence for any constant $c>0$.

The set of all 2iet sequences with an irrational slope $\alpha$ coincides with the set of all Sturmian sequences, see \cite{Lo2}.
The slope of the sequence $\uu$ with the vector of parameters $(\ell_0, \ell_1, \rho)^\top$ equals $ \alpha = \frac{\ell_1}{\ell_0+\ell_1}$.
The language of a Sturmian sequence depends only on the slope $\alpha$ and does not depend on the intercept $\delta$.  The following lemma summarizes  the relation between the slope and intercept of Sturmian sequences and the assigned vectors of parameters.

\begin{lem} \label{le:parameters_Sturmian}
Each lower Sturmian sequence ${\bf s}_{\alpha, \delta}$ is generated by the  transformation  $T$ exchanging  the intervals  $I_0 = [0, 1-\alpha)$ and  $I_1= [ 1-\alpha, 1)$; each  upper Sturmian sequence  ${\bf s'}_{\alpha, \delta}$ is generated by the  transformation  $T$ exchanging  the intervals  $I_0 = (0, 1-\alpha]$ and  $I_1= ( 1-\alpha, 1]$.

Moreover,

A) if $\delta \in (0,1)$, then
\[
 \vec{v}\left(  {\bf s'}_{\alpha, \delta} \right) = \vec{v}\left( {\bf s}_{\alpha, \delta}  \right) .
\]

B)  if $\delta = 0$, then
\[
\vec{v}\left( {\bf s'}_{\alpha, 0} \right) = \left( 1- \alpha, \alpha, 1  \right)
\quad \text{ and } \quad
\vec{v}\left( {\bf s}_{\alpha, 0}  \right) = \left( 1-\alpha, \alpha, 0 \right).
\]
\end{lem}

Among all Sturmian sequences with a  fixed irrational slope $\alpha=\frac{\ell_1}{\ell_0+\ell_1}$,   the sequence with the vector of parameters
$\vec{v}(\uu)= (\ell_0, \ell_1, \ell_1)^\top$  plays a special role.
Such a sequence is called a \textit{characteristic Sturmian sequence} and it is usually denoted by  ${\bf c}_\alpha$.
A Sturmian sequence $\uu \in \{0,1\}^\N$ is characteristic if both sequences $0\uu$ and $1\uu$ are Sturmian. Equivalently, a Sturmian sequence $\uu$ is characteristic, if every  prefix  of $\uu$ is  left special.

\section{Sturmian morphisms}

A morphism $\psi$ is a \emph{Sturmian morphism} if $\psi(\uu)$ is a Sturmian sequence for any Sturmian sequence $\uu$.
The set of Sturmian morphisms together with composition forms the so-called \emph{monoid of Sturm}, or \emph{Sturmian monoid}, which is in \cite{Lo2} denoted by ${\it St}$.
The monoid was described in \cite{MiPa_Rauzy}. Here we work  with a submonoid  of  ${\it St}$ generated by the following Sturmian morphisms:
\begin{equation}\label{listElementary}
  G: \begin{cases} 0 \to 0 \\ 1 \to 01 \end{cases} \quad
  \widetilde{G}: \begin{cases} 0 \to 0 \\ 1 \to 10 \end{cases} \quad
  D: \begin{cases} 0 \to 10 \\ 1 \to 1 \end{cases}
   \widetilde{D}: \begin{cases} 0 \to 01 \\ 1 \to 1 \end{cases} \quad .
\end{equation}
The submonoid $\mathcal{M} = \langle G,   \widetilde{G},  D, \widetilde{D} \rangle $ is also called the \emph{special Sturmian monoid}.
Any lower (resp. upper) mechanical sequence is mapped by a morphism from   $\mathcal{M}$  to a lower (resp. upper) mechanical sequence.
The monoid $\mathcal{M}$ does not contain the Sturmian morphism $E: 0\mapsto 1, 1\mapsto 0$. This morphism  maps  a lower (resp. upper) mechanical sequence to an upper (resp. lower) mechanical sequence.
Extending the generating set of $\mathcal{M}$ by the morphism $E$ gives already the generating set of the whole Sturmian monoid~${\it St}$.

Note that  $DE=EG$,  $\widetilde{D}E=E\widetilde{G}$  and  $E^2 = \mathrm{id}$.
Hence  ${\it St} = E \mathcal{M} \cup  \mathcal{M}  = \mathcal{M} E \cup  \mathcal{M} $.
In particular, $\psi^2 \in \mathcal{M}$ for each Sturmian morphism $\psi$. See also \Cref{rem:E} below.

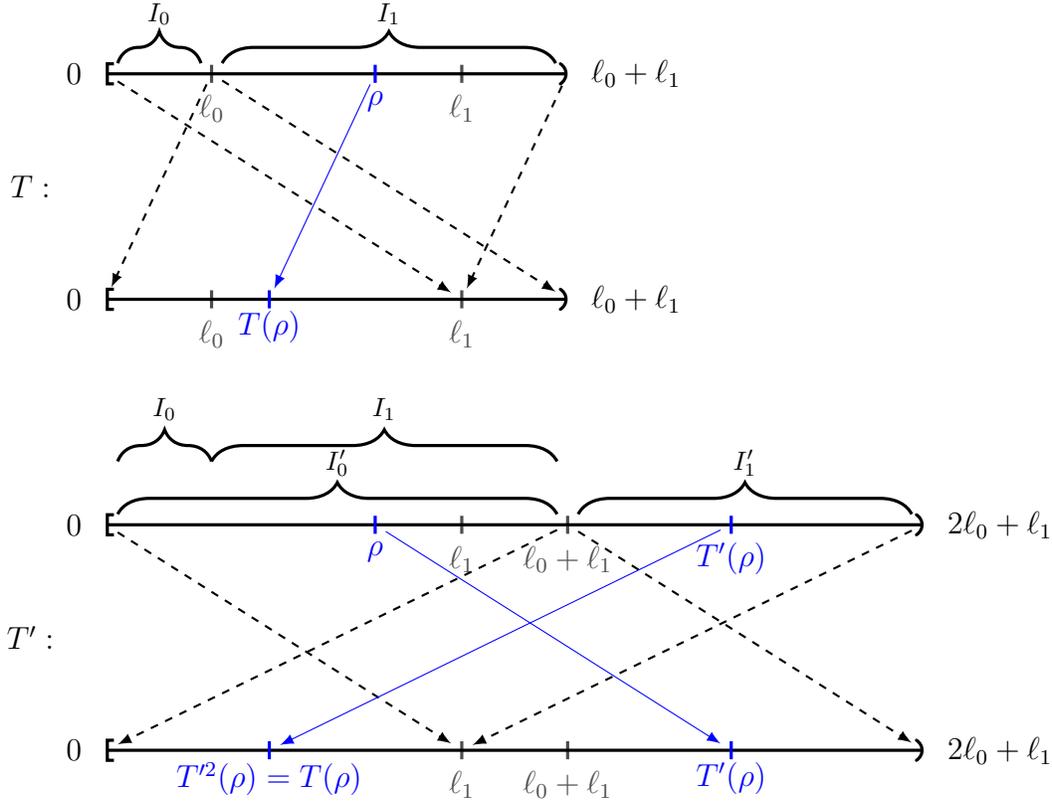
\begin{figure}
\newcommand*\ietalpha{11} 
\newcommand*\ietalphalabel{$\ell_0$} 
\newcommand*\ietbeta{37} 
\newcommand*\ietbetalabel{$\ell_1$} 
\newcommand*\ietX{28}

\begin{tikzpicture}[scale=1,radius=2cm,>=latex]

\newcommand\colorIntHeight{0.15}

\newcommand*\ietsdist{3}
\newcommand*\ietpad{0.2}
\newcommand*\ietepsstep{1.5}

\setcounter{ietstep}{1}
\setcounter{ietstepP}{0}
\setcounter{ietstepPP}{-1}

\tikzset{disco/.style={color=black,opacity=0.7}}
\tikzset{discoAl/.style={style=disco}}
\tikzset{discoBe/.style={style=disco}}
\tikzset{discoX/.style={color=blue}}
\tikzset{hlpo/.style={very thick}}
\tikzset{hlpoA/.style={yscale=1.2}}
\tikzset{tmap/.style={ultra thick,}}
\tikzset{bmap/.style={thick,dashed}}
\tikzset{xmap/.style={style=discoX},thick}
\tikzset{tmapG/.style={thick,dashed,color=gray}}

 \begin{scope}[scale=1,x=0.008\textwidth]

\node () at (103,0){};

\pgfmathsetmacro\newX{\ietX}

 \oneint[yshift={(1-\number\value{ietstep})*\ietsdist cm}]{\ietalpha}{\ietalphalabel}{\ietbeta}{\ietbetalabel}{\arabic{ietstep}}{%
  \ietpoint[pos=0.8,below]{discoX}{\newX}{$\rho$}{pointX\arabic{ietstep}%
   }
   \expandafter\draw [decorate,decoration={brace,amplitude=12pt,raise=4pt},yshift=0pt] (start\arabic{ietstep}) -- node[sloped,above=14pt] {\footnotesize $I_0$} (alpha\arabic{ietstep});
  \expandafter\draw [decorate,decoration={brace,amplitude=12pt,raise=4pt},yshift=0pt] (alpha\arabic{ietstep}) -- node[sloped,above=14pt] {\footnotesize $I_1$} (end\arabic{ietstep});

 }{$\ell_0+\ell_1$}

 \ietstepcounters

 \pgfmathsetmacro\newX{(\newX-\ietalpha)}

 \oneint[yshift={(1-\number\value{ietstep})*\ietsdist cm}]{\ietalpha}{\ietalphalabel}{\ietbeta}{\ietbetalabel}{\arabic{ietstep}}{%
  \ietpoint[pos=0.5,below]{discoX}{\newX}{$T(\rho)$}{pointX\arabic{ietstep}%
  }
 }{$\ell_0+\ell_1$}

\draw[->,bmap] (start\arabic{ietstepP}) to[]  node[pos=0.55,left,inner xsep=5] {} (beta\arabic{ietstep});
\draw[->,bmap] (alpha\arabic{ietstepP}) to[]  node[pos=0.55,left,inner xsep=5] {} (end\arabic{ietstep});

\draw[->,bmap] (alpha\arabic{ietstepP}) to[]  node[pos=0.55,left,inner xsep=5] {} (start\arabic{ietstep});
\draw[->,bmap] (end\arabic{ietstepP}) to[]  node[pos=0.55,left,inner xsep=5] {} (beta\arabic{ietstep});
\draw[->,xmap] (pointX\arabic{ietstepP}) to[]  node[pos=0.55,left,inner xsep=5] {} (pointX\arabic{ietstep});

\node[xshift={-1cm}] at ($(start\arabic{ietstepP})!0.5!(start\arabic{ietstep})$) {$T:$};

\pgfmathsetmacro\ietalphaOLD{\ietalpha}

\pgfmathsetmacro\ietalpha{\ietalpha+\ietbeta}
\renewcommand*\ietalphalabel{$\ell_0+\ell_1$} 
\pgfmathsetmacro\ietbeta{\ietbeta}
\renewcommand*\ietbetalabel{$\ell_1$} 

\ietstepcounters

\pgfmathsetmacro\newX{\ietX}
\pgfmathsetmacro\newTX{(\newX+\ietbeta)}
\pgfmathsetmacro\newTTX{(\newTX-\ietalpha)}

 \oneint[yshift={(1-\number\value{ietstep})*\ietsdist cm}]{\ietalpha}{\ietalphalabel}{\ietbeta}{\ietbetalabel}{\arabic{ietstep}}{%
  \ietpoint[pos=0.8,below]{discoX}{\newX}{$\rho$}{pointX\arabic{ietstep}%
   }
   \ietpoint[pos=0.8,below]{discoX}{\newTX}{$T'(\rho)$}{pointTX\arabic{ietstep}%
   }
   \expandafter\draw [decorate,decoration={brace,amplitude=12pt,raise=4pt},yshift=0pt] (start\arabic{ietstep}) -- node[sloped,above=14pt] {\footnotesize $I_0'$} (alpha\arabic{ietstep});
  \expandafter\draw [decorate,decoration={brace,amplitude=12pt,raise=4pt},yshift=0pt] (alpha\arabic{ietstep}) -- node[sloped,above=14pt] {\footnotesize $I_1'$} (end\arabic{ietstep});

  \begin{scope}[decoration={raise=24pt}]
  \draw [decorate,decoration={brace,amplitude=12pt}] (start\arabic{ietstep}) -- node[sloped,above=35pt] {\footnotesize $I_0$} (\ietalphaOLD,0);
  \draw [decorate,decoration={brace,amplitude=12pt}] (\ietalphaOLD,0) -- node[sloped,above=35pt] {\footnotesize $I_1$} (alpha\arabic{ietstep});
  \end{scope}

 }{$2\ell_0+\ell_1$}

 \ietstepcounters

 \pgfmathsetmacro\newX{(\newX+\ietbeta)}

 \oneint[yshift={(1-\number\value{ietstep})*\ietsdist cm}]{\ietalpha}{\ietalphalabel}{\ietbeta}{\ietbetalabel}{\arabic{ietstep}}{%
  \ietpoint[pos=0.5,below]{discoX}{\newX}{$T'(\rho)$}{pointX\arabic{ietstep}%
  }
  \ietpoint[pos=0.5,below]{discoX}{\newTTX}{$T'^2(\rho) = T(\rho)$}{pointTX\arabic{ietstep}%
   }
 }{$2\ell_0+\ell_1$}

\draw[->,bmap] (start\arabic{ietstepP}) to[]  node[pos=0.35,left,inner xsep=5] {} (beta\arabic{ietstep});
\draw[->,bmap] (alpha\arabic{ietstepP}) to[]  node[pos=0.35,left,inner xsep=5] {} (end\arabic{ietstep});

\draw[->,bmap] (alpha\arabic{ietstepP}) to[]  node[pos=0.35,left,inner xsep=5] {} (start\arabic{ietstep});
\draw[->,bmap] (end\arabic{ietstepP}) to[]  node[pos=0.35,left,inner xsep=5] {} (beta\arabic{ietstep});
\draw[->,xmap] (pointX\arabic{ietstepP}) to[]  node[pos=0.55,left,inner xsep=5] {} (pointX\arabic{ietstep});
\draw[->,xmap] (pointTX\arabic{ietstepP}) to[]  node[pos=0.55,left,inner xsep=5] {} (pointTX\arabic{ietstep});

\node[xshift={-1cm}] at ($(start\arabic{ietstepP})!0.5!(start\arabic{ietstep})$) {$T':$};

 \end{scope}
 \end{tikzpicture}
 \caption{Illustration of the relation of parameters of $\uu$ and $G(\uu)$: the upper part of the figure depicts 2iet transformation $T$ with the partition $[0,\ell_0)$ and $[\ell_0,\ell_0+\ell_1)$ which produces the sequence $\uu$.
 The lower part contains 2iet transformation $T'$ with the partition $I_0'=[0,\ell_0+\ell_1)$ and $I_1'=[\ell_0+\ell_1,2\ell_0+\ell_1)$ producing the sequence $G(\uu)$.
 }
 \label{obrazek}
 \end{figure}

The relation between the parameters of a Sturmian sequence and its image under a Sturmian morphism can be found in \cite[Lemmas 2.2.17 and 2.2.18]{Lo2}:
\begin{equation}\label{Lot}
G({\bf s}_{\alpha, \delta}) = {\bf s}_{\frac{\alpha}{1+\alpha}, \frac{\delta}{1+\alpha}}, \quad   \widetilde{G}({\bf s}_{\alpha, \delta}) = {\bf s}_{\frac{\alpha}{1+\alpha}, \frac{\alpha+ \delta}{1+\alpha}} \quad \text{and} \quad   E({\bf s}_{\alpha, \delta}) = {\bf s'}_{1-\alpha, 1-\delta}.
\end{equation}
The next lemma rephrases these  identities.
\Cref{obrazek} illustrates the action of the morphism $G$ on a Sturmian sequence and provides a geometrical proof of the first identity in \eqref{Lot}. Taking $x \in I_0$ in \Cref{obrazek} for $T$ and $T'$, we see that $T(x) = T'(x)$, while if $x \in I_1$, we have $T(x) = T'^2(x)$, $T'(x) \in I'_1$ and $T'(x) \not \in I_0 \cup I_1$.
 It follows that if $\uu$ is produced by the transformation $T$, then $G(\uu)$ is produced by the transformation $T'$.

\begin{lem}\label{image} Let $\uu$ be a lower (resp. upper) mechanical  sequence with parameters   $\ell_0 >0$, $\ell_1 >0$  and $\rho \in [0, \ell_0+\ell_1)$  (resp. $\rho \in (0, \ell_0+\ell_1]$). The lower  (resp.  upper) mechanical  sequence
\begin{itemize}
\item $G(\uu)$ has the parameters $\ell_0+\ell_1$, $\ell_1$ and $\rho $;

\item  $\widetilde{G}(\uu)$ has the parameters  $\ell_0+\ell_1$, $\ell_1$ and  $\rho +\ell_1$;

\item $D(\uu)$ has the parameters $\ell_0$, $\ell_0+\ell_1$ and  $\rho +\ell_0$;

\item  $\widetilde{D}(\uu)$ has the parameters $\ell_0$, $\ell_0+\ell_1$ and  $\rho $.
\end{itemize}
\end{lem}

\section{A representation of the special Sturmian monoid}\label{FR}

A \emph{representation} of the monoid $\mathcal{M}$ (over $\R$) is a monoid homomorphism $\mu: \mathcal{M} \to \R^{n \times n}$.
A representation is \emph{faithful} if $\mu$ is injective.

A traditional representation of the monoid $\mathcal{M}$ (and in general, of a monoid of morphisms) assigns to a morphism  $\varphi \in \mathcal{M}$ its \emph{incidence matrix} $M_\varphi \in \mathbb{N}^{2\times2}$ defined by
\[
(M_\varphi) _{i,j} = |\varphi(j)|_i \quad \text{ for } i,j\in  \{0,1\},
\]
i.e., the entry of $M$ at the position $(i,j)$  equals the number of occurrences of the letter $i$ in the word $\varphi(j)$.
It is easy to verify that $M_{\varphi\circ \psi} = M_{\varphi}M_{\psi}$ for every pair $\varphi, \psi  \in \mathcal{M}$, hence $\varphi \mapsto M_\varphi$ is a representation of the special Sturmian monoid $\mathcal{M}$.

However, the matrices assigned to $\widetilde{G}$ and $G$ coincide.
The same is true for the matrices assigned to $\widetilde{D}$ and $D$.
In other words, the representation $\varphi \mapsto M_\varphi$ is not injective and hence it is not faithful.

Let us recall that the special Sturmian monoid $\mathcal{M}$ is not free: for any $k \in \N$ we have
\begin{equation}\label{eq:relations}
     GD^k\widetilde{G} = \widetilde{G}\widetilde{D}^kG \quad \text{and} \quad    DG^k\widetilde{D} = \widetilde{D}\widetilde{G}^kD.
\end{equation}
Theorem 2.3.14 in~\cite{Lo2} says that \eqref{eq:relations} is the presentation of the monoid $\mathcal{M}$, i.e., no other non-trivial independent relation can be found.

\Cref{image} is the core of the faithful representation of $\mathcal{M}$ we construct below.  It  says that the action of a Sturmian morphism on Sturmian sequences corresponds to an action of a linear operator on  vectors from $\mathbb{R}^3$.
To specify a faithful representation of the monoid $\mathcal{M}$
we assign to each element of $\{G,\widetilde{G},D, \widetilde{D}\}$ one matrix:
\begin{equation}\label{matrices}
R_{\widetilde{G}} =\left( \begin{array}{lll}
1&1&0\\
0&1&0\\
0&1&1
\end{array}\right), \
R_G =\left( \begin{array}{lll}
1&1&0\\
0&1&0\\
0&0&1
\end{array}\right), \
R_{\widetilde{D}} =\left( \begin{array}{lll}
1&0&0\\
1&1&0\\
0&0&1
\end{array}\right),  \ \text{ and } \
R_D =\left( \begin{array}{lll}
1&0&0\\
1&1&0\\
1&0&1
\end{array}\right).
\end{equation}
These matrices preserve the presentation of the special Sturmian monoid \eqref{eq:relations}.

\begin{claim}\label{aswell}
If  $k \in \N$, then
\begin{equation}\label{eq:MatrixRelations}
R_{\widetilde{G}} R_{\widetilde{D}}^k R_G =   R_{G} R_D^k R_{\widetilde{G}} \quad \text{ and } \quad R_{\widetilde{D}} {R_{\widetilde{G}}^k} R_D  = R_{D} R_G^k {R_{\widetilde{D}}}.
\end{equation}
\end{claim}
\begin{proof} For each $k\in \N$ we have
$$
R_{\widetilde{G}}^k =\left( \begin{array}{lll}
1&k&0\\
0&1&0\\
0&k &1
\end{array}\right), \
R_{G}^k =\left( \begin{array}{lll}
1&k&0\\
0&1&0\\
0&0&1
\end{array}\right), \
R_{\widetilde{D}}^k =\left( \begin{array}{lll}
1&0&0\\
k&1&0\\
0&0&1
\end{array}\right),  \ \text{ and } \
R_D^k =\left( \begin{array}{lll}
1&0&0\\
k&1&0\\
k&0&1
\end{array}\right).
$$

Hence,
\[
R_{\widetilde{G}} R_{\widetilde{D}}^k R_G =
\begin{pmatrix}
  k+1&k+2&0\\
k&k+1&0\\
k&k+1&1
  \end{pmatrix}  = R_{G} R_D^k R_{\widetilde{G}} \quad \text{and} \quad R_{\widetilde{D}} {R_{\widetilde{G}}^k} R_D =
  \begin{pmatrix}
   k+1&k&0\\
k+2&k+1&0\\
k+1&k &1
   \end{pmatrix} = R_{D} R_G^k {R_{\widetilde{D}}}.\qedhere
\]
\end{proof}

We can now assign a matrix to any element of $\mathcal{M}$.

\begin{defi} \label{def:R}
Let $\mathcal{R}: \mathcal{M} \mapsto \R^{3 \times 3}$ be defined
for $\psi \in \mathcal{M}$ by
\[
\mathcal{R}(\psi) = R_{\varphi_1}R_{\varphi_1}\cdots R_{\varphi_n},
\]
where $\psi= \varphi_1 \circ \varphi_2\circ \cdots \circ \varphi_n$ and $\varphi_i \in \{G,\widetilde{G},D, \widetilde{D}\}$ for every $i = 1,2, \ldots, n$.
\end{defi}

Let us note that the definition is correct. It does not depend on the decomposition of $\psi$ into the elements of $\{G,\widetilde{G},D, \widetilde{D}\}$
since by \Cref{aswell} the relations  \eqref{eq:relations}  of the presentation of $\mathcal{M}$ are preserved in the monoid
$\left \langle R_G,R_{\widetilde{G}},R_D, R_{\widetilde{D}}  \right \rangle = \mathcal{R}(\mathcal{M})$.

Our choice of the four matrices in \eqref{matrices} enables the following matrix reformulation of \Cref{image}.

\begin{claim}\label{obrazyVektoru}
Let $\varphi \in \mathcal{M}$ and $ \vec{v}(\uu)$ be the vector of parameters of a Sturmian sequence  $ \uu$.
The vector $\mathcal{R}(\varphi) \vec{v}(\uu)$ is the vector of parameters of the Sturmian sequence $\varphi(\uu)$.
\end{claim}
\begin{proof}
Let $\varphi = \varphi_1 \circ \varphi_2\circ \cdots \circ \varphi_n$ and $\varphi_i \in \{G,\widetilde{G},D, \widetilde{D}\}$ for every $i = 1,2, \ldots, n$.
We proceed by induction on $n$.
For $n = 0$, the claim trivially holds.
Set $\varphi_2 \circ \varphi_3\circ \cdots \circ \varphi_{n} = \varphi'$  and
assume that $\mathcal{R}(\varphi') \vec{v}(\uu)$ is the vector of parameters of $\varphi'(\uu)$.
Since $\varphi_1 \in \{G,\widetilde{G},D, \widetilde{D}\}$, then $R_{\varphi_1} \left( \mathcal{R}(\varphi') \vec{v}(\uu)  \right) $ is the vector of parameters of $\varphi_1 \left( \varphi'(\uu) \right) $ by \Cref{image}.
The proof is finished by noticing $R_{\varphi_1} \mathcal{R}(\varphi')  = \mathcal{R}(\varphi)$.
\end{proof}

As already mentioned, the traditionally used representation of the Sturmian monoid which maps a morphism to its incidence matrix is not faithful.
Nevertheless, this representation has a strong connection to the representation $\mathcal{R}$.
By definition of the matrices in \eqref{matrices}, $\mathcal{R}$ maps the morphism $\psi$ to the matrix $ \mathcal{R}(\psi)$ of the form
\begin{equation}\label{submatrix}\mathcal{R}(\psi) = \left( \begin{array}{lll}
m_{00}&m_{01}&0\\
m_{10}&m_{11}&0\\
E&F&1
\end{array}\right), \  \text{ where }  \left( \begin{array}{ll}
m_{00}&m_{01}\\
m_{10}&m_{11}
\end{array}\right) = M_\psi\ \ \text{ and }  \ \  E,F \in \mathbb{N}.
\end{equation}
In other words, the incidence matrix of $\psi$ is the top left submatrix of $\mathcal{R}(\psi)$.

\begin{proposition}\label{DefRepr}
 The mapping $\mathcal{R}$ is a faithful representation of the monoid $\mathcal{M}$.
\end{proposition}
\begin{proof}
\Cref{aswell} and  the presentation of $\mathcal{M}$ by \eqref{eq:relations} imply that $\mathcal{R}$ is a representation of $\mathcal{M}$.
To show injectivity of $\mathcal{R}$, we assume $\mathcal{R}(\psi) = \mathcal{R}(\varphi)$, for $\psi, \varphi \in \mathcal{M}$.
Let $\uu$ be a lower Sturmian sequence.
By \Cref{obrazyVektoru} and \Cref{image} both morphisms $\psi$ and $\varphi$ map $\uu$ to lower Sturmian sequences with the same parameters.
That is, $\varphi(\uu) = \psi(\uu)$.
By the same argument,  $\varphi(\uu) = \psi(\uu)$  for each upper Sturmian sequence $\uu$, too. In other words,  the images by the morphisms $\psi$  and $\varphi$ coincide for each element of their domain.  We conclude that $\varphi = \psi$. Hence, $\mathcal{R}$ is injective.
\end{proof}

We continue by listing several straightforward properties of the representation $\mathcal{R}$ of the special Sturmian monoid.
Let us remind that a \emph{convex cone} $C$ in $\mathbb{R}^n$ is a non-empty subset of $\mathbb{R}^n$ such that
$C\cap (-C) =\{0\}$ and $\mu x  + \nu y \in C$ for each $x, y \in C$ and $\mu, \nu \geq 0$.
For example, the set
\[
\R_{\geq 0}^n = \{x \in \mathbb{R}^n: \text{all components of $x$ are non-negative} \}
\]
is a convex cone.

\begin{lem}\label{list}
Let $\psi \in \mathcal{M}$.
\begin{enumerate}
\item  The convex cone $C_1: =\{(x,y,z)^\top \in \mathbb{R}^3: 0 \leq x, \ 0\leq y, \  0\leq z\leq x+y\}$  is  invariant under multiplication by $\mathcal{R}(\psi)$, i.e.,
$\mathcal{R}(\psi) (C_1) \subset C_1$.
\medskip
\item The convex cone $C_2 := \{(x,y,z)^\top \in \mathbb{R}^3 :  0 \leq x,\ 0\geq y, \ y\leq z\leq x\}$  is  invariant under multiplication by the inverse matrix  of  $\mathcal{R}(\psi)$, i.e.,
$\bigl(\mathcal{R}(\psi)\bigr)^{-1} (C_2) \subset C_2$.

\medskip

\item  The number $1$ is an eigenvalue of $\mathcal{R}(\psi)$ and $(0,0,1)^\top$ is  its corresponding eigenvector.
In particular, the convex cone $C_3 := \{(0,0,z)^\top \in \mathbb{R}^3 :  0\leq z\}$  is  invariant under multiplication by $\mathcal{R}(\psi)$.
\end{enumerate}
\end{lem}
 \begin{proof}
The validity of  all items for the matrices $R_G, R_{\widetilde{G}}, R_D, R_{\widetilde{D}}$ can be verified directly.
Since $\mathcal{R}(\psi)$ belongs to the monoid generated by these four matrices,  $\mathcal{R}(\psi)$ satisfies these properties as well.
\end{proof}

In order to provide some applications of the representation $\mathcal{R}$, we recall three  properties of a linear mapping preserving a closed convex cone. The first item mentioned in the following proposition is a consequence of the Brouwer's theorem, see for instance \cite{DuSch}; the second and the third item are a consequence of the Perron-Frobenius theory, see for instance \cite{Fiedler}.
Let us note that the assumption $A\left(\R_{\geq 0}^n\right) \subset \R_{\geq 0}^n$ in second item below is equivalent to the property that all entries of $A$ are non-negative.

\begin{prop}\label{PerFro} Let $A\in \mathbb{R}^{n\times n}$  and $C $ be a closed  convex cone in $\mathbb{R}^n$.

\begin{enumerate}[(1)]

\item If  $AC \subset C$, then at least one eigenvector of $A$ belongs to the convex cone  $C$. \label{it:PF1}

\item If  $A \left(\R_{\geq 0}^n\right) \subset \R_{\geq 0}^n$, then  the spectral radius of $A$ is an eigenvalue corresponding to an eigenvector from $\R_{\geq 0}^n$. \label{it:PF2}

\item If all entries of the matrix  $A^k$ are positive for some $k \in \mathbb{N}$, then
 the spectral radius $r_A$ of $A$ is a dominant simple eigenvalue of $A$, i.e. all other eigenvalues are of modulus strictly smaller than $r_A$. The corresponding eigenvector to $r_A$ has all entries positive and an  eigenvector corresponding to  any other eigenvalue cannot have all entries non-negative. \label{it:PF3}
 \end{enumerate}
\end{prop}

\begin{coro}\label{eigenvectors}
Let $\psi \in \mathcal{M}$ be a primitive morphism.
The matrix $\mathcal{R}(\psi)$  has eigenvalues $\Lambda$, $1$ and $\frac{1}{\Lambda}$, where $\Lambda >1$ is a quadratic unit.   An eigenvector corresponding to  $\Lambda$ can be found in the form  $(x,y,z)^\top \in \bigl(\mathbb{Q}(\Lambda)\bigr)^3$ with $x>0$, $y>0$ and $z\geq 0$.
No other eigenvalue has an eigenvector with the first two components positive.
\end{coro}

\begin{proof}
The matrix  $\mathcal{R}(\psi)$ has the form described in  \eqref{submatrix}.
The vector $(0,0,1)^\top$ is its eigenvector corresponding to the eigenvalue $1$.
Other eigenvalues of $\mathcal{R}(\psi)$ are eigenvalues of the matrix  $M_\psi$ as well.
Since all entries of  $M_\psi$ are non-negative integers and $\det M_\psi = 1$,  the two eigenvalues $\lambda_1$ and  $\lambda_2$ of  $M_\psi$  are  roots of the polynomial $X^2 - pX+1$ with $p=Tr(M_\psi) >0$.
Thus  $\lambda_1$, and $ \lambda_2$ are  algebraic integers and $\lambda_1 \lambda_2=1$.
Primitivity of $\psi$ implies that for some  integer  $k\geq 1$  all entries of   $M_\psi^k$  are positive.
By \Cref{it:PF3} of \Cref{PerFro}, the spectral radius of $M_\psi$ is a simple eigenvalue.
It follows, without loss of generality, $\Lambda :=\lambda_1 > 1 > \lambda_2>0$ and $\Lambda$ is the dominant eigenvalue of $\mathcal{R}(\psi)$.
Since $\lambda_2$ is an algebraic integer lying in the interval $(0,1)$, $\lambda_2$ cannot be rational.
Consequently, $\lambda_2$ and $\lambda_1$ are quadratic irrational numbers.
By \Cref{it:PF2} of \Cref{PerFro},  $\mathcal{R}(\psi)$ has a non-negative eigenvector, say  $(x,y,z)^\top$, corresponding  to $\Lambda$.
Since  $\bigl(\mathcal{R}(\psi)-\Lambda I\bigr) (x,y,z)^\top = (0,0,0)^\top$, the entries $x,y,z$ can be chosen to belong into $\mathbb{Q}(\Lambda)$.
The vector $(x,y)^\top$ is an eigenvector of $M_\psi$ corresponding to $\Lambda$ and by \Cref{it:PF3} of \Cref{PerFro}, $x> 0, y>0$.

If $(x',y',z')$ is an eigenvector of $\mathcal{R}(\psi)$ corresponding to $\lambda_2$, then $(x',y')$ is an eigenvector of $M_\psi$ to the same  eigenvalue.
It follows from \Cref{it:PF3} of \Cref{PerFro} that $x'$ and $y'$ are numbers with  opposite signs.

\end{proof}

\begin{remark} \label{rem:E} The Sturmian morphism $E$ exchanging the letters $0\leftrightarrow 1$ can be associated with a matrix $R_E\in \Z^{3\times 3}$ such that $R_E \vec{v}(\uu)$ is a vector of parameters of the Sturmian sequence $E(\uu)$.
Using \eqref{Lot}, we see that
$$R_E =  \left( \begin{array}{llc}
0&1&0\\
1&0&0\\
1&1&-1
\end{array}\right).$$
Therefore,  it is possible to extend the representation $\mathcal{R}$ of the special Sturmian monoid $\mathcal{M}$ into a representation of the Sturmian monoid ${\it St}$.
Since we focus on the question  when  a Sturmian sequence is invariant under a primitive morphism, we can restrict ourselves to $\mathcal{M}$ only: indeed, if $\uu$ is fixed by a primitive morphism  $\psi\in {\it St}$,  then  $\uu$ is fixed by the  morphism $ \psi^2\in \mathcal{M}$.
Using this restriction, we avoid having negative entries in the matrices representing the morphisms and we can exploit the Perron-Frobenius theorem.
\end{remark}

\section{A submonoid of $Sl(\mathbb{Z},3)$ defined by three convex cones}\label{submonoidSl}

Let us recall the notation
$Sl(\mathbb{N},3) = \{ R \in \mathbb{N}^{3\times 3} : \det R = 1\}$ and  $Sl(\mathbb{Z},3) = \{ R \in \mathbb{Z}^{3\times 3} : \det R = 1\}$.   It is well known that the group  $Sl(\mathbb{Z},3)$ is finitely generated. Rivat  in \cite{Fogg} showed   that  the monoid $Sl(\mathbb{N},3)$ is not finitely generated. Clearly,  our representation $\mathcal{R}$  maps the special Sturmian monoid $\mathcal{M}$  into $Sl(\mathbb{N},3)$. In particular, $\mathcal{R}(\mathcal{M})$ is a submonoid of $Sl(\mathbb{N},3)$ that is finitely generated by the four matrices $R_G, R_{\widetilde{G}}, R_D, R_{\widetilde{D}}$.    Lemma \ref{list} suggests to study   a new  submonoid of $Sl(\mathbb{Z},3)  $.  We set
\begin{equation}\label{newMonoid}\mathcal{E} = \{ R \in Sl(\mathbb{Z},3) : RC_1 \subset C_1, \   R^{-1}C_2\subset C_2, \  RC_3 \subset C_3\}, \end{equation}
where
\[
\begin{split}
C_1 &:= \{(x,y,z)^\top \in \mathbb{R}^3 \colon  0 \leq x, 0\leq y,  0\leq z\leq x+y\}, \\
C_2 &:= \{(x,y,z)^\top \in \mathbb{R}^3 \colon  0 \leq x,0\geq y, y\leq z\leq x\}, \text{ and } \\
C_3 &:= \{(0,0,z)^\top \in \mathbb{R}^3 \colon  0\leq z\}
\end{split}
\]
(in accordance with the notation in \Cref{list}).

In this section we show that the representation  $\mathcal{R}$ is in fact an isomorphism between  $\mathcal{M}$  and $\mathcal{E}$.
By \Cref{list},   $\mathcal{R}(\mathcal{M})\subset \mathcal{E} $. In the remaining part of the section we prove
$\mathcal{R}(\mathcal{M})=\mathcal{E} $.
To achieve this goal, we  characterize the elements of $\mathcal{E}$ by inequalities.
We first give a technical lemma that shall be used to reduce the number of inequalities.

\begin{lem}\label{1enough} If $A,B,C,D,E,F \in \mathbb{N}$ satisfy  \ \
$  E < A+C, \ \  F<B+D\,$  and $ AD-BC = 1$, then
$$ -A < AF -BE \leq B \quad \Longleftrightarrow\quad   -C \leq CF-DE < D\,.$$
\end{lem}

\begin{proof}
$(\Longrightarrow)$   If $C=0$,  then $AD-BC =1$  implies  $A=D=1$  and $E< A+C$  gives  $E < 1$, i.e. $E=0$.  The desired inequality  $-C \leq CF-DE < D$ has  now the form $0\leq 0< 1$ and it is obviously satisfied.

 If  $C>0$, then the restated assumption $-A+1 \leq  AF -BE \leq B$ implies
$$-AC+C \leq  AFC -BEC \leq BC \,.$$
We continue by substituting $BC =AD-1$
$$-AC +C \leq AFC - E(AD-1) \leq AD -1\  \Longrightarrow  \
-AC +C -E\leq A(CF - ED) \leq AD -1-E.  $$
As $0\leq E\leq A+C-1$, we obtain
$$-AC -A+1\leq A(CF - ED) <  AD  \ \  \Longrightarrow  \ \   -C    -\tfrac{A-1}{A}\leq CF -ED < D,
$$
where we used $A>0$ which is a simple consequence of $AD-BC = 1$.
As $A,C,D,E,F \in \mathbb{N}$  and $\tfrac{A-1}{A} \in [0,1)$, the last inequality gives the desired $-C \leq CF-DE < D$.

\bigskip

$(\Longleftarrow)$    \ \ Let us note that the relations   $  E < A+C, \ \  F<B+D\,$  and $ AD-BC = 1\,$   in the assumption of the lemma are invariant  under exchange of the parameters $A \leftrightarrow D$,
$C \leftrightarrow B$ and $E \leftrightarrow F$. Applying the exchange to the implication
$ -A < AF -BE \leq B \quad \Longrightarrow\quad   -C \leq CF-DE < D$ \  which is already proven, we obtain the converse.
\end{proof}

\begin{lem}\label{inequalities}
If $R \in Sl(\mathbb{Z},3)$, then $R \in \mathcal{E}$ if and only if  there exist $A,B,C,D,E,F \in \mathbb{N}$ such that
$$R=\begin{pmatrix}
A&B&0\\
C&D&0\\
E&F&1
\end{pmatrix}, \quad \text{where}  \
AD-BC = 1 \quad  \text{and}$$
\begin{equation}\label{straight}
 E < A+C, \ \  F<B+D\,,
\end{equation}
\begin{equation}\label{lower}
 -C \leq CF-DE < D\,.
\end{equation}
\end{lem}

\begin{proof} $(\Longrightarrow)$ \ \ Let   $R=\begin{pmatrix}
A&B&G\\
C&D&H\\
E&F&J
\end{pmatrix}\in \mathcal{E}$, where $A,B,\ldots , H, J\in \mathbb{Z}$.
Since $RC_3\subset C_3$,  we have $R  \left( \begin{array}{l}
0\\
0\\
1
\end{array}\right)=\left( \begin{array}{l}
G\\
H\\
J
\end{array}\right) \in C_3$. Therefore,     $G=H=0$ and $J\geq 0$.   The determinant of $R$ is 1 and can be now computed as $\det R = (AD-BC)J=1$.  Since $ A,B,C,D, J$ are  integers  and $J\geq 0$,  we have $AD-BC=1=J$.
\medskip

By the assumption $RC_1\subset C_1$ we have
$$\left( \begin{array}{lll}
A&B&0\\
C&D&0\\
E&F&1
\end{array}\right)  \left( \begin{array}{l}
x\\
y\\
z
\end{array}\right)=\left( \begin{array}{l}
Ax+By\\
Cx+Dy\\
Ex+Fy+z
\end{array}\right) \in C_1   \text{ \ \  for each }\left( \begin{array}{l}
x\\
y\\
z
\end{array}\right)\in C_1.$$

The inequalities $Ax+By\geq 0$, $Cx+Dy\geq 0$  and $Ex+Fy+z\geq 0$  for all choices $x,y\geq 0$  and $z=0$ imply  $A,B, C, D,E,F\geq 0$, i.e., $R \in Sl(\mathbb{N},3)$.  Moreover,
$ Ex+Fy+z \leq Ax+By +Cx+Dy$ holds  for any $z$ satisfying $0 \leq z \leq x + y$. In particular, if $z = x+y$, then
$$
Ex+Fy +x+y\leq Ax+By+Cx+Dy \ \ \text{gives} \ \ 0\leq (A+C-E-1)x+ (B+D-F-1)y.
$$
Therefore, $A+C-E-1\geq 0$ and $B+D-F-1\geq 0$.

\medskip

Since $R^{-1}C_2\subset C_2$ and $R^{-1} = \begin{pmatrix}
D&-B&0\\
-C&A&0\\
FC-ED&BE-FA&1
\end{pmatrix} $, we have
$$
R^{-1} \begin{pmatrix}
x\\
y\\
z
\end{pmatrix}
 = \begin{pmatrix}
Dx-By\\
-Cx+Ay\\
(FC-ED)x+(BE-FA)y+z
\end{pmatrix} \in C_2
$$
for any choice $y\leq 0\leq x$ and $y\leq z\leq x$.
In particular, the  third coordinate satisfies two inequalities:
\begin{enumerate}[i)]
\item $-Cx+Ay \leq (FC-ED)x + (BE-FA)y +z $, and \label{it:ine1}
\item $(FC-ED)x + (BE-FA)y +z \leq Dx-By $. \label{it:ine2}
\end{enumerate}

Setting $ y=z=0$ and $x \geq 0$, inequality \ref{it:ine1} implies $FC-ED+C \geq 0$ .
Similarly, putting $0\leq z =x $ and $y=0$ we obtain $FC-ED-D+1\leq 0$ from \ref{it:ine2}.

$(\Longleftarrow)$: We have to verify 3 conditions in the definition of $\mathcal{E}$ in \eqref{newMonoid}.   Let us start with checking $R^{-1}C_2 \subset C_2$. We  consider $(x,y,z)^\top \in C_2$, i.e., $x\geq 0$, $y\leq  0$, $y\leq z\leq x$.  Then
$$
R^{-1} \begin{pmatrix}
x\\
y\\
z
\end{pmatrix}
 = \begin{pmatrix}
Dx-By\\
-Cx+Ay\\
(FC-ED)x+(BE-FA)y+z
\end{pmatrix} = :\begin{pmatrix}
\widetilde{x}\\
\widetilde{y}\\
\widetilde{z}
\end{pmatrix} .  $$
Obviously,  the first component $ \widetilde{x} = Dx-By \geq 0$, as $x\geq 0\geq y$ and $B,D\geq 0$. By the same reason,  the second component $ \widetilde{y} = -Cx+Ay\leq 0$.  To deduce inequalities required  for the third component $\widetilde{z}$ we use \eqref{lower} and by   \Cref{1enough}  also the inequality   $ -A < AF -BE \leq B$. Hence
$$ 0\geq \underbrace{(FC-ED-D+1)}_{\leq 0}x - \underbrace{(-B-BE+FA)}_{\leq 0}y = \widetilde{z} - \widetilde{x}  +\underbrace{x -z}_{\geq 0}  \geq  \widetilde{z} - \widetilde{x}.
$$
Analogously,
$$ 0\leq \underbrace{(FC-ED+C)}_{\geq 0}x + \underbrace{(-A+1 +BE-FA)}_{\leq 0}y = \widetilde{z} - \widetilde{y}  +\underbrace{y-z}_{\leq 0}  \leq  \widetilde{z} - \widetilde{y}.
$$
We can conclude  that $(\widetilde{x},\widetilde{y},\widetilde{z})^\top \in C_2$.

To  verify that $RC_1\subset C_1$ and $RC_3\subset C_3$ is straightforward and we omit it. \end{proof}

Before we show that $\mathcal{E}$ equals the image of the special Sturmian monoid $\mathcal{M}$ under our representation~$\mathcal{R}$, we recall a property of
$Sl(\mathbb{N},2) = \{ M \in \mathbb{N}^{2\times 2} : \det M = 1\}$.

\begin{claim}\label{size2} If $M=
\begin{pmatrix}
A&B\\
C&D
\end{pmatrix}
\in Sl(\mathbb{N},2)$ and $M\neq I$, then either ($A\geq C$  and $B\geq D$) or  ($A\leq C$  and $B\leq D$).
\end{claim}
\begin{proof}
Let us note that $A>0$ and $D>0$. Indeed, if $AD=0$, then $1=\det M= AD-CB \leq 0$ --- a contradiction.

Assume that the claim does not hold true, i.e., ($A < C$  or $B < D$) and  ($A > C$ or $B > D$).

First discuss the case $A < C$.
Necessarily,  $B>D$. Since $1=\det M = AD-BC \leq (C-1)D - (D+1)C= -D-C \leq 0$ --- a contradiction.

Now  let us discuss the case  $B < D$. Hence $A>  C$.
It follows that $1=\det M = AD-BC \geq (C+1)D - (D-1)C= D+C \geq 1$.
Therefore, we can replace inequalities by equalities. As $D>0$ and $  C\geq 0$, the last inequality gives   $C=0$ and $D=1$.
Consequently, $A=C+1 = 1$ and $B=D-1 =0$.  In other words $M=I$ --- a contradiction.
\end{proof}

\begin{thm}\label{rovnost}
The monoid  $\mathcal{E}$  defined by \eqref{newMonoid} coincides with
$\mathcal{R}(\mathcal{M}) = \left \langle R_G, R_{\widetilde{G}}, R_D, R_{\widetilde{D}} \right \rangle$.
\end{thm}

\begin{proof}
The inclusion  $\mathcal{R}(\mathcal{M})\subset \mathcal{E} $ follows from \Cref{list}.

Let us show $\mathcal{E} \subset \mathcal{R}(\mathcal{M})$.
Let $R = (R_{ij})\in \mathcal{E}$ have the form given by \Cref{inequalities}, i.e.,
\[
R=\begin{pmatrix}
A&B&0\\
C&D&0\\
E&F&1
\end{pmatrix}
\]
and all entries of $R$ are non-negative.

First  assume that $C=0$.   As $1=\det R =AD-BC$, we have $A=D=1$. Inequalities  \eqref{straight} say  $E=0$ and $0\leq F\leq B$. Therefore,  $R$ has the form
$$R= \begin{pmatrix}
1&B&0\\
0&1&0\\
0&F&1
\end{pmatrix} = {\underbrace{\begin{pmatrix}
1&1&0\\
0&1&0\\
0&0&1
\end{pmatrix}}_{R_G}}^{B-F}
{\underbrace{\begin{pmatrix}
1&1&0\\
0&1&0\\
0&1&1
\end{pmatrix}}_{R_{\widetilde G}}}^F \in \mathcal{R}(\mathcal{M})\,.$$
Now   assume that $B =0$.   Analogously to the previous case,  $\det R =1$ implies $A=D=1$. Inequalities  \eqref{straight} give  $0\leq E\leq C$,  $F=0$  and
$$R= \begin{pmatrix}
1&0&0\\
C&1&0\\
E&0&1
\end{pmatrix} = {\underbrace{\begin{pmatrix}
1&0&0\\
1&1&0\\
0&0&1
\end{pmatrix}}_{R_{\widetilde D}}}^{C-E}
{\underbrace{\begin{pmatrix}
1&0&0\\
1&1&0\\
1&0&1
\end{pmatrix}}_{R_D}}^E \in \mathcal{R}(\mathcal{M})\,.$$

We proceed by induction on $R_{11} + R_{21} = A+C$. Note that $\det R=1$ forces  $A,D >0$.
Hence the case $A+C = 1$ implies $A=1,C=0$ and it is treated above.
Further, in the induction step, it suffices to discuss only the situation  $A,B,C,D \geq 1$.
Using \Cref{size2}, we split our discussion into the following two cases:

\begin{enumerate}[(I)]
\item $A\geq C\geq 1$ and $B\geq D\geq 1$, \label{it:pr_case_1}
\item $1\leq A\leq  C$ and $1\leq B\leq D$. \label{it:pr_case_2}
\end{enumerate}
\medskip

\noindent  {\bf Case \ref{it:pr_case_1}} $A\geq C\geq 1$ and $B\geq D\geq 1$:

First we show  by contradiction that entries of $R$ satisfy
\begin{equation}\label{subcases}
\left(C\leq E \ \text{ and } \ D\leq F \right)
\quad \text{ or } \quad
\left (E< A \ \text{ and } \ F< B\right).
\end{equation}

Assuming that \eqref{subcases} does not hold and using inequalities of Case~\ref{it:pr_case_1} we obtain
\[
\left ( C>E \ \text{ and } \ F\geq B\geq D  \geq 1 \right)
\quad \text { or } \quad
\left ( B\geq D>F \ \text{ and } \ E\geq A\geq C \geq 1 \right).
\]
In the first case, we have $C \geq E+1$ and $D \leq F$.
Combining these inequalities with  inequality \eqref{lower} of \Cref{inequalities} we obtain
\[
D \stackrel{\eqref{lower}}{>} CF -DE \geq  (E+1)F - FE = F\geq D,
\]
which is a contradiction.

In the second case, we have $F \leq B - 1$ and $E \geq A$.
By~\Cref{inequalities}, the assumptions of~\Cref{1enough} are satisfied, and inequality \eqref{lower}
can be restated as $ -A < AF -BE \leq B$, hence
\[-A < AF -BE.
\]
Together with $F \leq B - 1$ and $E \geq A$ we conclude
\[
-A < AF - BE \leq A(B-1) - BA = - A,
\]
which is a contradiction.

\medskip
\noindent {\bf Subcase (I.a) } $C\leq E$ and $D\leq F$:

We observe that the matrix $R'$ defined  via
$$R= \begin{pmatrix}
A&B&0\\
C&D&0\\
E&F&1
\end{pmatrix} =
\underbrace{\begin{pmatrix}
1&1&0\\
0&1&0\\
0&1&1
\end{pmatrix}}_{R_{\widetilde{G}}}
\underbrace{\begin{pmatrix}
A-C&B-D&0\\
C&D&0\\
E-C&F-D&1
\end{pmatrix}}_{=:R'}$$ is in $Sl(\mathbb{N},3)$, as all entries of $R'$ are non-negative integers and $\det R' = 1$.

We continue by proving $R' \in \mathcal{E} $ using \Cref{inequalities} for $R'$.
Inequalities \eqref{straight} for the matrix $R'$ read $E-C < A-C+C$  and $F-D < B-D +D$.
These are equivalent to $E<A+C$ and $F< B+D$, which are satisfied due to $R \in \mathcal{E}$ and \Cref{inequalities}.
Inequalities  \eqref{lower} for $R'$ say $-C\leq C(F-D) - (E-C)D < D$.
They are equivalent to $-C\leq CF-ED < D$, again satisfied by $R \in \mathcal{E}$ and \Cref{inequalities}.
Therefore $R' \in \mathcal{E}$.

Let us recall that  we proceed by induction on  $R_{11}+R_{21}$.
Since  $R'_{11}+R'_{21} = A<  A+C = R_{11}+R_{21}$, we know by the induction hypothesis that  $R' \in \left \langle R_G, R_{\widetilde{G}}, R_D, R_{\widetilde{D}} \right \rangle$.
The relation $R=R_{\widetilde{G}}R'$ implies $R\in \left\langle R_G, R_{\widetilde{G}}, R_D, R_{\widetilde{D}} \right\rangle$  as well.

\medskip
\noindent {\bf Subcase (I.b) } $E< A$ and $F< B$:

We define the matrix $R'$ by
$$R= \begin{pmatrix}
A&B&0\\
C&D&0\\
E&F&1
\end{pmatrix} =
\underbrace{\begin{pmatrix}
1&1&0\\
0&1&0\\
0&0&1
\end{pmatrix}}_{R_G}
\underbrace{\begin{pmatrix}
A-C&B-D&0\\
C&D&0\\
E&F&1
\end{pmatrix}}_{=:R'}.$$
Obviously, $ R' \in Sl(\mathbb{N},3)$.
We use again \Cref{inequalities} to show that $R' \in \mathcal{E} $.
Inequality \eqref{straight}  for $R'$ states $E< A-C+C$ and $F< B-D +D$.
These are satisfied by the specification Subcase I.b.
Inequalities \eqref{lower} for $R'$ and $R$ coincide and they follow from $R\in \mathcal{E}$ and \Cref{inequalities}.
As $R'_{11}+R'_{21}< R_{11}+R_{21}$ and $R=R_GR'$, we conclude that $R\in \langle R_G, R_{\widetilde{G}}, R_D, R_{\widetilde{D}} \rangle$.

\medskip

\noindent  {\bf Case \ref{it:pr_case_2}}  $1\leq A\leq C $ and $1\leq B\leq  D$:

We transform this case to the Case~\ref{it:pr_case_1}.
We use the permutation matrix  $P=\begin{pmatrix}
0&1&0\\
1&0&0\\
0&0&1
\end{pmatrix}$ and the following 3 facts:
\begin{itemize}
\item $P\langle R_G, R_{\widetilde{G}}, R_D, R_{\widetilde{D}} \rangle P = \langle R_G, R_{\widetilde{G}}, R_D, R_{\widetilde{D}} \rangle$.
\\ \emph{Proof:} The equality follows from the equalities $PR_{\widetilde{G}}P = R_D$,  $PR_GP = R_{\widetilde{D}}$  and $P^2 =I_3$.

\item $P\mathcal{E}P =\mathcal{E}$.
\\ \emph{Proof:} The cones $C_1$, $C_2$ and $C_3$ satisfy $PC_1 = C_1$, $PC_2 =-C_2$ and $PC_3 = C_3$.
Let $R\in \mathcal{E}$.
We have $PRP \in Sl(\mathbb{Z},3)$ and  $  (PRP)^{-1}C_2 = PR^{-1}P C_2 = -PR^{-1}C_2 \subset - PC_2 = C_2$.
Analogously, $PRPC_i \subset C_i $ for $i$ being $1$ or $3$.
It means that  $PRP \in \mathcal{E}$, or equivalently $\mathcal{E} \supset P\mathcal{E}P$.
Multiplying the last inclusion by $P$ from the right and left, we obtain $P\mathcal{E}P \supset P^2\mathcal{E}P^2 = \mathcal{E}$.

\item If $R\in \mathcal{E}$  belongs to Case~\ref{it:pr_case_2}, then $PRP$ belongs to Case~\ref{it:pr_case_1}.
\end{itemize}
 Proof of the theorem is now complete. \end{proof}

\section{The representation  $\mathcal{R}$ and fixed points of Sturmian morphisms}\label{app}
In this section we apply the faithful representation $\mathcal{R}$ to study the parameters of fixed points of primitive morphisms.

\subsection{Parameters of a fixed point  of Sturmian morphisms }
In article \cite{Peng} Peng and Tan solve the question to determine the parameters of a fixed point $\uu$ of a given primitive Sturmian morphism $\psi$.
To find the slope of $\uu$, i.e. the frequency of the letter $1$ in $\uu$,  one can use  the incidence matrix of $M_\psi$.
For a primitive morphism,  it is well-known that the positive components of the eigenvector corresponding to the dominant eigenvalue are proportional to the frequencies of letters.
Hence, the only  non-trivial question is to determine the intercept of $\uu$. \
Theorems 2.3 and 3.2 of \cite{Peng} answer this question using the structure of the words $\psi(01)$ and $\psi(10)$.
The faithful representation $\mathcal{R}(\psi)$ provides  a simple algebraic method to determine the parameters of $\uu$.

\begin{prop}\label{fixed}
 Let  $\psi\in  \mathcal{M}$  be   a primitive morphism and $\uu$ be a Sturmian sequence with the vector of parameters $\vec{v}(\uu)$.
 The sequence $\uu$  is  fixed by $\psi$ if and only if  $\vec{v}(\uu)$  is an eigenvector to the dominant eigenvalue of $\mathcal{R}(\psi)$.
\end{prop}

\begin{proof}
Let the sequence $\uu$ with parameters $\vec{v}(\uu) = (\ell_0, \ell_1,\rho)^\top$ be fixed by $\psi $.
Since the Sturmian sequences $\uu$ and $\psi(\uu)$ coincide, the vectors of their parameters are collinear, i.e. there exists $\Lambda > 0$ such that  $\vec{v}(\psi(\uu)) = \Lambda \vec{v}(\uu)$.
By \Cref{obrazyVektoru},
\[
\Lambda \vec{v}(\uu) = \vec{v}(\psi(\uu)) =\mathcal{R}(\psi)\vec{v}(\uu).
\]
By \Cref{eigenvectors}, an eigenvector $ (\ell_0, \ell_1,\rho)^\top$ with positive components $\ell_0$ and $ \ell_1$ corresponds to the dominant eigenvalue of $\mathcal{R}(\psi)$.

To prove the converse, assume that $\vec{v}(\uu)$  is an eigenvector corresponding to the dominant eigenvalue $\Lambda$ of $\mathcal{R}(\psi)$.
By \Cref{obrazyVektoru},  the sequence $\psi(\uu)$ has  the vector of parameters $\vec{v}(\psi(\uu)) =\mathcal{R}(\psi)\vec{v}(\uu) = \Lambda \vec{v}(\uu)$, i.e., the vectors of parameters of $\uu$ and $\psi(\uu)$ are the same up to a scalar factor.
Recall that any morphism of $\mathcal{M}$  maps  a lower  (resp. upper) mechanical sequence  to a lower  (resp. upper) mechanical sequence.
Two  lower  (resp. upper)  mechanical sequences with the same vector of parameters up to a scalar factor coincide.
Consequently, $\psi(\uu) = \uu$.
\end{proof}

\subsection{Pairs ${\bf s}_{\alpha, \delta}$ and  ${\bf s'}_{\alpha, \delta }$ fixed by (possibly distinct) morphisms}

Dekking \cite{Dekking} studies for which values of the slope and the intercept are both sequences ${\bf s}_{\alpha, \delta}$ and  ${\bf s'}_{\alpha, \delta}$ fixed by primitive morphisms.
His result can be also proven by applying the representation  $\mathcal{R}$.
Recall that the definition of  ${\bf s}_{\alpha, \delta}$ and  ${\bf s'}_{\alpha, \delta}$ immediately  implies that  ${\bf s}_{\alpha, \delta}$ and  ${\bf s'}_{\alpha, \delta}$  are either identical or they differ  at most on two  neighbouring positions.

First we state several simple claims on invariant subspaces of  matrices from $\mathcal{R}(\mathcal{M})$.
A subspace $V \subset \R^3$ is an \emph{invariant subspace} of a matrix $R$ if  $RV \subset V$.
By inspecting the  behaviour of the matrices assigned to the generators of $\mathcal{M}$ we obtain the following properties.

\begin{lem}\label{subspaces} Let  $\psi \in \mathcal{M}$.
\begin{enumerate}

\item If $\psi \in \langle G, D \rangle $,   then   $\mathcal{R}(\psi)$ has an invariant subspace
 $\{(x,y,z)^\top :\ y=z\} \subset \mathbb{R}^3$.

\item  If $\psi \in \langle \widetilde{G}, \widetilde{D} \rangle $,   then    $\mathcal{R}(\psi)$ has an invariant subspace   $\{(x,y,z)^\top :\ x=z\} \subset \mathbb{R}^3$.

\item If    $\psi \in \langle \widetilde{G}, D \rangle $,   then   $\mathcal{R}(\psi)$ has an invariant subspace   $\{(x,y,z)^\top :\ z=x+y\} \subset \mathbb{R}^3$.

\item  If  $\psi \in \langle G, \widetilde{D} \rangle $,   then  $\mathcal{R}(\psi)$ has an invariant subspace   $\{(x,y,z)^\top :\ z=0\} \subset \mathbb{R}^3$.

\end{enumerate}
\end{lem}

The last lemma is transformed directly into the next statement in terms of parameters of a Sturmian sequence.

\begin{lem}\label{subspaces2}
Let  a Sturmian sequence  $\uu$ with parametrs $\ell_0, \ell_1, \rho$ be fixed by a primitive  morphism $\psi  \in \mathcal{M}$.
\begin{enumerate}

\item      If  $\psi \in \langle \widetilde{G}, D \rangle $, then     $\rho =\ell_0+\ell_1$.

\item If   $\psi \in \langle G, D \rangle $, then   $\rho =\ell_1$.

\item   If    $\psi \in \langle \widetilde{G}, \widetilde{D} \rangle $, then    $\rho =\ell_0$.

\item If     $\psi \in \langle G, \widetilde{D} \rangle $, then    $\rho =0$.

\end{enumerate}

\end{lem}

\begin{proof}  Item (1): Assume $\psi \in \langle \widetilde{G}, D \rangle $.   By Proposition \ref{fixed}, the vector $(\ell_0, \ell_1, \rho)^\top$ is the eigenvector of $\mathcal{R}(\psi)\in \langle R_{\widetilde{D}}, R_G\rangle$ corresponding to the dominant eigenvalue.
By \Cref{subspaces}, the plane $P=\{(x,y,z)^\top :\ z=x+y\} \subset \mathbb{R}^3$ is invariant under multiplication by $R_G$ and by $R_{\widetilde{D}}$, therefore $\mathcal{R}(\psi)$ has two eigenvectors in the plane. As the eigenvector $(0,0,1)^\top$  of  $\mathcal{R}(\psi)$ corresponding to $1$ does not belong to  $P$,  the eigenvector corresponding to the dominant eigenvalue $(\ell_0, \ell_1, \rho)^\top$ must be in $P$.
It implies $\rho = \ell_0+\ell_1$, and hence  the Sturmian sequence $\uu$ is an upper Sturmian sequence coding the two interval exchange transformation with the domain  $(0, \ell_0+\ell_1]$.

Proofs of the other Items are analogous.
\end{proof}

Now we state Dekking's (see \cite[Theorems 2 and 3]{Dekking}) result and provide its alternative proof.

\begin{prop}
Let $\alpha \in (0,1)$, $\alpha$ irrational, and $\delta \in [0,1)$.
Assume that both sequences  ${\bf s}_{\alpha, \delta}$ and  ${\bf s'}_{\alpha, \delta}$ are fixed by primitive morphisms  and ${\bf s}_{\alpha, \delta}\neq {\bf s'}_{\alpha, \delta}$.  Then
either

\begin {enumerate}
\item $\delta = 1-\alpha$, in which case  ${\bf s}_{\alpha, \delta }$ and  ${\bf s'}_{\alpha, \delta}$  are distinct fixed points  of the same  primitive morphism $\psi  \in \langle \widetilde{G}, \widetilde{D} \rangle$; or
\item $\delta =0$, in which case  ${\bf s}_{\alpha, \delta}$  is fixed by a morphism $\psi \in \langle G, \widetilde{D} \rangle$ and  ${\bf s'}_{\alpha, \delta}$  is fixed by a morphism $\eta \in \langle \widetilde{G}, D \rangle$.
Moreover, if $\psi = \varphi_1\circ \varphi_2 \circ \cdots \circ \varphi_n$ with  $\varphi_i\in  \{G, \widetilde{D}\}$, then
$$ \eta = \xi_1\circ\xi_2\circ \cdots \circ \xi_n, \ \ \text{where} \ \  \xi_i=\begin{cases} \widetilde{G} & \text{if } \ \varphi_i =G, \\ D & \text{if } \ \varphi_i=\widetilde{D}, \end{cases} \quad  \text{for every } \  i=1, \ldots, n.
$$
\end{enumerate}
\end{prop}

\begin{proof}
Since ${\bf s}_{\alpha, \delta}$ is fixed by a primitive morphism, there exists a primitive morphism $\psi\in \mathcal{M}$ which fixes ${\bf s}_{\alpha, \delta}$ (see \Cref{rem:E}).


Firstly, we assume that $\delta \in (0,1)$.
Due to \Cref{le:parameters_Sturmian},  $\vec{v}({\bf s}_{\alpha, \delta}) = \vec{v}({\bf s'}_{\alpha, \delta}) = (1-\alpha, \alpha, \delta)^\top$.
By \Cref{fixed},  ${\bf s'}_{\alpha, \delta}$  is fixed by  $\psi $ as well. As  ${\bf s}_{\alpha, \delta}\neq {\bf s'}_{\alpha, \delta}$, the primitive  morphism  $\psi$ has two fixed points.  In particular,  $\psi(0)$ has a prefix $0$ and  $\psi(1)$ has a prefix $1$.

The form of morphisms  $ G,D, \widetilde{G}, \widetilde{D}  $ implies that the starting letters of $\psi(0)$ and $\psi(1)$   coincide whenever the  morphism  $G$ or $D$ occurs in the composition of $\psi$. Therefore, $\psi \in \langle \widetilde{G}, \widetilde{D} \rangle$.
By  Lemma \ref{subspaces2},  $\delta = 1-\alpha$.

\medskip

Secondly, we assume that $\delta = 0$.  Due to \Cref{le:parameters_Sturmian},  $\vec{v}({\bf s}_{\alpha, \delta}) =  (1-\alpha, \alpha,0)^\top$ and   $ \vec{v}({\bf s'}_{\alpha, \delta}) = (1-\alpha, \alpha, 1)^\top$.

Let us recall  that   $M_\psi$ is the product of incidence matrices of the elementary morphisms $M_D = M_{\widetilde{D}}$ and $M_G = M_{\widetilde{G}}$.
Let us write  $M_\psi = M_{\varphi_1}M_{\varphi_2} \cdots M_{\varphi_n}$,  where ${\varphi_i} \in \{ G, \widetilde{D}\}$, and  define $\varphi = \varphi_1\circ \varphi_2\circ \cdots \circ \varphi_n$. Obviously, $M_\psi = M_\varphi$ and $\mathcal{R}(\varphi) =  R_{\varphi_1}R_{\varphi_2}\cdots R_{\varphi_n}$.   Since the third row of both the matrices $R_G$  and $R_{\widetilde{D}}$ equals $(0,0,1)$, the third row of the matrix $\mathcal{R}(\varphi)$ is $(0,0,1)$ as well. Hence $\mathcal{R}(\psi) = \mathcal{R}(\varphi)$, c.f. \eqref{submatrix}.  As the representation is faithful, $\psi = \varphi \in \langle G, \widetilde{D}\rangle$.

Let $\eta $  be  the morphism created in Item (2). The incidence matrices $M_\eta$ and $M_\psi$ coincide, since $M_G=M_{\widetilde{G}}$ and  $M_D=M_{\widetilde{D}}$. Obviously,  $
 (1-\alpha, \alpha)^\top$ is the positive eigenvector of $M_\eta$ as well. By Lemma \ref{subspaces}, the positive eigenvector of  $\mathcal{R}(\eta)$ belongs to the plane formed by the vectors $(x,y,z)^\top$ for which $z=x+y$. Therefore,  the positive eigenvector of  $\mathcal{R}(\eta)$  equals $(1-\alpha, \alpha, 1)^\top$. By Proposition \ref{fixed},  ${\bf s'}_{\alpha, \delta}$ is fixed by $\eta$.
\end{proof}

\subsection{Sturmian sequences fixed by  a morphism}
Yasutomi in \cite{Ya99} gave a necessary and sufficient condition under which a Sturmian sequence is invariant under a primitive morphism. The same result is proved in \cite{BaMaPe} and \cite{BEIR}.

Using the faithful representation $\mathcal{R}$ of the Sturmian monoid we provide here a simple proof that the condition is necessary.
Due to the relation $E({\bf s}_{\alpha, \delta}) = {\bf s'}_{1-\alpha, 1-\delta}$ mentioned in  \eqref{Lot}, it is enough to characterize  lower Sturmian sequences fixed by a primitive morphism.

\begin{proposition} \label{Yasutomi_necessary} Let $\alpha \in (0,1)$ be irrational  and $\delta \in [0,1)$.  If   ${\bf s}_{\alpha, \delta}$ is  fixed by a primitive morphism, then
\begin{enumerate}
\item
 $\alpha$ and $ \delta$ belong to the same quadratic field, say $\mathbb{Q}(\sqrt{m})$;
 \item $\overline{\alpha} \notin (0,1) $ and  $\min \{ \overline{\alpha}, 1-\overline{\alpha}\} \leq  \overline{\delta} \leq    \max  \{ \overline{\alpha}, 1-\overline{\alpha} \}$, where the mapping $x\mapsto \overline{x}$ is the non-trivial field automorphism on  $\mathbb{Q}(\sqrt{m})$ induced by $\sqrt{m} \mapsto - \sqrt{m}$.
\end{enumerate}
\end{proposition}
\begin{proof}   Let $\psi$ be a primitive morphism fixing ${\bf s}_{\alpha, \delta}$. The vector $\vec{v} = (1-\alpha, \alpha, \delta)^\top$ is a vector of parameters of the Sturmian sequence  ${\bf s}_{\alpha, \delta}$.   By  Proposition \ref{fixed},   $\vec{v}$ is an eigenvector to the dominant eigenvalue $\Lambda$ of the matrix $\mathcal{R}(\psi)$.  Due to Corollary \ref{eigenvectors}:
\begin{itemize}

\item  $\Lambda$ is a quadratic number, i.e. $\mathbb{Q}(\Lambda) = \mathbb{Q}(\sqrt{m})$ for some $m \in \mathbb{N}$.

\item   $c\vec{v} = c(1-\alpha, \alpha, \delta)^\top  \in  \bigl(\mathbb{Q}(\Lambda)\bigr)^3$ for some positive $c$.  Consequently, $c =c(1-\alpha) + c \alpha \in  \mathbb{Q}(\Lambda)$. It implies  $(1-\alpha, \alpha, \delta)^\top  \in  \bigl(\mathbb{Q}(\Lambda)\bigr)^3$ and thus $\delta$ and  $\alpha$ are in the same quadratic field $\mathbb{Q}(\sqrt{m})$.

\end{itemize}
Let us apply  the field automorphism to $\mathcal{R}(\psi) (1-\alpha, \alpha, \delta)^\top = \Lambda (1-\alpha, \alpha, \delta)^\top$. Since the entries of $\mathcal{R}(\psi)$ are rational, $\mathcal{R}(\psi)$  is invariant under the automorphism and thus  we get that
 $(1-\overline{\alpha}, \overline{\alpha}, \overline{\delta})^\top$ is an eigenvector to the eigenvalue $\overline{\Lambda}$. The third eigenvector of  $\mathcal{R}(\psi)$ is $(0,0,1)^\top$.

 By Item 2) of Lemma \ref{list} and Item 1) of  Proposition \ref{PerFro},  one eigenvector of $\bigl(\mathcal{R}(\psi)\bigr)^{-1}$  belongs to the cone $C_2 := \{(x,y,z)^\top \in \mathbb{R}^3 :  0 \leq x,\ 0\geq y, \ y\leq z\leq x\}$.
 Now we use the fact that if  a non-singular matrix   has an eigenvector  $\vec{d}$ to $\lambda$, then  $\vec{d}$ is an eigenvector of its inverse matrix  to the eigenvalue $\tfrac{1}{\lambda}$. Hence   the eigenvector   $ (1-\overline{\alpha}, \overline{\alpha}, \overline{\delta})$ or $ -(1-\overline{\alpha}, \overline{\alpha}, \overline{\delta})$  belongs to $C_2$.
\medskip

If $ (1-\overline{\alpha}, \overline{\alpha}, \overline{\delta})^\top \in C_2$, then $\overline{\alpha} \leq \overline{\delta} \leq 1-\overline{\alpha}$.

If $(-1+\overline{\alpha}, -\overline{\alpha}, -\overline{\delta})^\top \in C_2$, then $-\overline{\alpha} \leq -\overline{\delta} \leq -1+\overline{\alpha}$, or equivalently, $1-  \overline{\alpha} \leq \overline{\delta} \leq \overline{\alpha}$.

Both cases confirm   Item (2) of the proposition. \end{proof}
Let us stress that Yasutomi also proved that  Items (1)
 and (2) of the previous proposition  are also sufficient for  ${\bf s}_{\alpha, \delta}$   to be fixed by a primitive morphism.

 \begin{remark}\label{kombinace} Yasutomi's characterization implies that
      if  ${\bf s}_{\alpha, \delta_1}$  and ${\bf s}_{\alpha, \delta_2}$ are fixed by primitive morphisms, then  ${\bf s}_{\alpha, \delta}$  is fixed by a primitive morphism for every $\delta = c\delta_1 + (1-c)\delta_2$, where $c \in (0,1)\cap\mathbb{Q}$. Indeed, $\overline{\delta} = c\overline{\delta_1} + (1-c)\overline{\delta_2}$  belongs to the interval  $\bigl[\min \{ \overline{\alpha}, 1-\overline{\alpha}\},   \max  \{ \overline{\alpha}, 1-\overline{\alpha} \}\bigr]$ since $\overline{\delta}$ is a convex combination of two numbers  from this interval.

 \end{remark}

\subsection{Conjugacy of Sturmian morphisms}

In this subsection we use the faithful representation of the special Sturmian monoid to deduce a known result on the number of morphisms which are conjugates of a given Sturmian morphism. First we  recall  some notions.

We say that a morphism $\varphi$
is a \emph{right conjugate} of a morphism $\psi$, or that $\psi$ is a \emph{left
conjugate} of $\varphi$,  noted $\psi\triangleright\varphi$, if there
exists $w \in \A^*$ such that
\begin{equation}
 w\psi(a) = \varphi(a)w, \quad \textrm{for every letter   } a \in \A. \label{FirstCond}
\end{equation}
If \eqref{FirstCond} is satisfied for $\psi = \varphi$ with a non-empty $w$, then the morphism $\psi$ is called \emph{cyclic}.
A fixed point  of a cyclic morphism is periodic,   hence  Sturmian morphisms are acyclic.  To any   acyclic morphism one may   assign a morphism $\varphi_R$, called the  \emph{rightmost conjugate of $\varphi$}, such
that the following two conditions hold:
\begin{enumerate}[\rm (i)]
\item $\varphi_R$ is a right conjugate of $\varphi$;
\item if $\xi$ is a right conjugate of $\varphi_R$, then $\xi=\varphi_R$.
\end{enumerate}
\begin{remark}\label{leftspecial}
Let us list some simple properties of the relation $\triangleright$.
\begin{enumerate}
    \item If $\psi\triangleright\varphi$, then $|\varphi(a)|_b = |\psi(a)|_b$ for every $a,b \in \mathcal{A}$. Hence the incidence matrices  of $\varphi$ and $\psi$ coincide, i.e.,  $M_{\varphi} = M_{\psi}$. In particular, $\psi$ is primitive if and only if $\varphi$ is primitive.

    \item  Let  $\psi: \mathcal{A}^*\mapsto \mathcal{A}^*$ be a  morphism and $\uu$ be an infinite sequence over $\mathcal{A}$ such that every letter of $\mathcal{A}$ occurs in $\uu$. If $\psi\triangleright\varphi$ and $\uu =\varphi(\uu) = \psi(\uu)$, then $\varphi = \psi$.

\item If an acyclic morphism $\varphi$ acts on binary alphabet $\{0,1\}$,      the last letters of $\varphi_R(0)$ and   of $\varphi_R(1)$  are distinct. Consequently,  if $0u$ and $1u$ occur in a fixed point of $\varphi_R$, then $0\varphi_R(u)$ and $1\varphi_R(u)$ occur in the fixed point as well.
\end{enumerate}
\end{remark}

The following result can be found in \cite[Proposition 2.3.21]{Lo2}.

\begin{theorem} \label{sturmian_conjugate}  If $M=\begin{pmatrix}
A&B\\
C&D
\end{pmatrix} \in Sl(\N,2)$, then $M$ is the incidence matrix of  $A+B+C+D-1$ mutually conjugate Sturmian morphisms.
\end{theorem}

\begin{proof} Fix  $S \in \{ 0, 1, \ldots, A+B+C+D-2\}$.  First we show that there exists a unique pair  $(E,F) \in \N^2$ such that  $E+F = S$ and  $R= \begin{pmatrix}
A&B&0\\
C&D&0\\
E&F&1
\end{pmatrix}  \in \mathcal{E}$.  By \Cref{inequalities,1enough} we look for $E \in \N$ and $F=S-E$ satisfying $-A < A(S-E) -BE \leq B$. Or equivalently,   $\frac{AS-B}{A+B}\leq E < \frac{AS+A}{A+B}$. Since the distance between the lower and the upper bounds on $E$ equals 1, exactly one integer $E$, namely $E = \lceil \frac{AS-B}{A+B}\rceil$,  satisfies $-A < A(S-E) -BE \leq B$.  Let us check that such $E$ and $F:=S-E$ satisfy  \eqref{straight}.
Using  the definition of $E$ and   $A\geq 1$ we obtain
$$  -1<\tfrac{-B}{A+B} \leq\tfrac{AS-B}{A+B}\leq E < \tfrac{AS-B}{A+B} + 1 \leq \tfrac{A(A+B+C+D-2)-B}{A+B} + 1 = A+C -\tfrac{A-1}{A+B}\leq A+C.
$$
As $E$ is an integer, the previous inequalities confirm that $0\leq E < A+C$.
To check \eqref{straight} for  $F$,  we use $F:=S-E$  and write
\begin{multline*}
    -1\leq\tfrac{-A}{A+B}\leq \tfrac{SB-A}{A+B}=S - \tfrac{AS-B}{A+B} -1< F\leq S - \tfrac{AS-B}{A+B} = \\
    = \tfrac{(S+1)B}{A+B}\leq  \tfrac{B(A+B+C+D-1)}{A+B} = B+D -\tfrac{1+B}{A+B}.
\end{multline*}
It confirms that $0\leq F < B+D$.
By \Cref{1enough},  the $E$ and $F$ satisfy \eqref{lower} as well.

Observe that if  $S \notin \{0,1,\ldots, A+B+C+D-2 $\},  then no values $E,F \in \mathbb{N}, S=E+F$,   satisfy inequalities $E<A+C $  and   $F<B+D $ required by  \Cref{inequalities}.     We have proved that any $M \in Sl(\N,2)$  occurs in the left upper corner of  a matrix from $\mathcal{E}$ exactly  $A+B+C+D -1$ times.
By \Cref{rovnost}  and the fact that the representation  $\mathcal{R}$ of the monoid $\mathcal{M}$ is faithful, there exist in  $\mathcal{M}$ exactly $A+B+C+D-1$ Sturmian morphisms  having the incidence matrix $ M$.      Let us denote these morphisms by $\varphi^{(i)}$ for $i=0, 1, \ldots, A+B+C+D-2$.  We need to show that these morphisms are mutually conjugate, or equivalently,   to show that  the rightmost conjugate  $\varphi_R^{(i)}$ does not depend on the index $i$.

If the matrix $M$ is primitive, then the frequencies of letters in a fixed point of  $\varphi^{(i)}$ form  an eigenvector to the dominant eigenvalue of $M$. Therefore, all fixed points of the  $A+B+C+D-1$ morphisms
have the same frequencies of letters, i.e., they have the same slope, say  $\alpha$. By Item (3) of  Remark   \ref{leftspecial}, any prefix of the fixed point of  $\varphi_R^{(i)}$ is a left special factor of the fixed point. In other words, the fixed point of      $\varphi_R^{(i)}$ is the characteristic sequence  ${\bf c}_\alpha$ for every index $i$.  By Item (2) of the same remark,   $\varphi_R^{(i)} = \varphi_R^{(0)}$ for every $i$.

If the matrix $M$ is not primitive, then either $B=0$ or $C=0$. Assume that $C=0$. It follows that the  morphism $\varphi^{(i)}$ is of  the form  $0\mapsto 0, 1\mapsto 0^i10^{C-i}$, for $i = 0, 1, \ldots, C$. It is obvious that they are mutually conjugate. The case $B=0$ is analogous.
\end{proof}

\section{The square root of fixed point of characteristic Sturmian morphisms}\label{ctverce}
Saari~\cite{Sa} showed that for every Sturmian sequence $\uu$ there exist 6 its
factors $w_1,\ldots,w_6$ such that \begin{equation}\label{squares} \uu = w_{i_1}^2 w_{i_2}^2 w_{i_3}^2
\ldots, \qquad \text{ where $i_k \in
\{1,\ldots,6\}$ for each } k\in\mathbb{N},\end{equation}  and moreover, for each $k\in \N$, the shortest square prefix of the sequence $w_{i_k}^2 w_{i_{k+1}}^2
w_{i_{k+2}}^2\dots$ is $w_{i_{k}}^2$.
This result served as inspiration for J. Peltomäki and M. Whiteland to introduce the square root $\sqrt{\uu}$ of the Sturmian sequence $\uu$ written in the form \eqref{squares} as $$ \sqrt{\uu} = w_{i_1}
w_{i_2} w_{i_3} \ldots$$
In \cite{PeWh},  they also proved the following theorem.
\begin{theorem}\label{pelto} If  $\uu$  is a Sturmian sequence  with the slope $\alpha$ and the intercept $\delta$, then
the sequence $\sqrt{\uu}$ is a Sturmian sequence with the  same slope $\alpha$ and the intercept
$\frac{1-\alpha+\delta}{2}$.
\end{theorem}
\begin{example}\label{priklad}
    Let $\varphi\in\M$ such that $\varphi = DG^2: 0\mapsto 10, 1 \mapsto 10101$.
    The fixed point of $\varphi$ can be written as concatenation of the squares  of these 6 factors:  $10$, $1$, $0110101$, $101$, $01$, $01101$.
    The beginning of this decomposition is as follows:
    \begin{align*}
        \uu &=
        10101101010110101011010110101011010101101010110101101010\dots\\
        &=\underbrace{10}_{w_1} \underbrace{10}_{w_1} \underbrace{1}_{w_2}
        \underbrace{1}_{w_2} \underbrace{01}_{w_3} \underbrace{01}_{w_3}
        \underbrace{0110101}_{w_4} \underbrace{0110101}_{w_4}
        \underbrace{10}_{w_1} \underbrace{10}_{w_1}
        \underbrace{101}_{w_5} \underbrace{101}_{w_5}
        \underbrace{01}_{w_3} \underbrace{01}_{w_3}
        \underbrace{10}_{w_1} \underbrace{10}_{w_1}\dots
    \end{align*}
    Hence, the square root begins with
    \begin{align*}
        \sqrt{\uu} &=
        101010110101101010110101011010101101011010101
        \dots\\
        &=\underbrace{10}_{w_1}
        \underbrace{1}_{w_2}
        \underbrace{01}_{w_3}
        \underbrace{0110101}_{w_4}
        \underbrace{10}_{w_1}
        \underbrace{101}_{w_5}
        \underbrace{01}_{w_3}
        \underbrace{10}_{w_1}
        \underbrace{101}_{w_5}
        \underbrace{0110101}_{w_4}
        \underbrace{0110101}_{w_4}
        \underbrace{10}_{w_1}
        \underbrace{101}_{w_5}
        \underbrace{01}_{w_3} \dots
    \end{align*}
\end{example}

\begin{theorem}\label{naseodmocniny2} Let $\uu\in\{0,1\}^\mathbb{N}$ be a  Sturmian sequence
    fixed by a primitive morphism $\varphi\in\M$. The square root $\sqrt{\uu}$ is fixed by
    a morphism $\psi\in\M$ which is a conjugate of one of the morphisms $\varphi,
    \varphi^2,    \varphi^3$ or $\varphi^4$.  \end{theorem}

\begin{proof} Let $\vec{v}(\uu)=
    (1-\alpha, \alpha, \delta)^\top$ be the vector of parameters of  $\uu$.
  According to \Cref{fixed},  the
    vector $\vec{v}(\uu)$ is an eigenvector of $\Rcal(\varphi)$ corresponding to the dominant eigenvalue. Set $P:= \begin{pmatrix} 1 & 0 & 0 \\ 0 & 1 & 0 \\
    \frac{1}{2} & 0 &  \frac{1}{2}\end{pmatrix}$.  Then by  Theorem \ref{pelto} we have $\vec{v}(\sqrt{\uu})=
    (1-\alpha, \alpha, \frac{1-\alpha+\delta}{2} )^\top = P\vec{v}(\uu)$.
    Clearly $P\vec{v}(\uu)$ is
    the eigenvector of the matrix $P\Rcal(\varphi^k)P^{-1}$
    for every $k \in \mathbb{N}$
    and by Corollary~\ref{eigenvectors}
    the eigenvector
    corresponds to the dominant eigenvalue.

For the moment, assume  \begin{equation}\label{cochceme}P\Rcal(\varphi^k)P^{-1}\in \mathcal{E} = \mathcal{R}(\mathcal{M})\quad \text{for some }\ k \in \{1,2,3,4\}\,.\end{equation}
    It means that  $\sqrt{\uu}$ is fixed by a morphism $\psi \in
    \mathcal{M}$, for which  $\mathcal{R}(\psi) = P\mathcal{R}(\varphi^k)P^{-1}$. Because of  the form of the matrix $P$, the left upper
    $(2\times 2)$-submatrices of the matrices  $\mathcal{R}(\psi)$ and $\mathcal{R}(\varphi^k)$ coincide and they are equal to the incidence matrices of the morphisms $\psi$ and $\varphi^k$. By Theorem \ref{sturmian_conjugate}, the morphism $\psi$ is conjugate to the morphism $\varphi^k$ as stated.

    Therefore, to complete the proof we have to verify  $\eqref{cochceme}$.
    Let $Q\in\N^{3\times 3}$ denote
    $Q = \Rcal(\varphi) = \left(\begin{smallmatrix}
        M & \vec{0} \\ e_1^{\top} & 1
    \end{smallmatrix}\right)$ where $M=M_\varphi$ is the incidence matrix
    of the morphism $\varphi$ and $e_1^{\top} = (E,F)$.
    Then $\Rcal(\varphi^k) = \Rcal(\varphi)^k = Q^k$.
    We prove $\eqref{cochceme}$ in 4 steps.

    \textit{i)}: It holds that $Q^k = \left(\begin{smallmatrix}
        M^k & \vec{0} \\ e_k^{\top} & 1
    \end{smallmatrix}\right)$ where $e_k^{\top}
    = e_1^{\top}\sum_{i=0}^{k-1} M^i$. Indeed, from the relation
    $Q^{k+1}=Q^k Q$ we get the recurrence relation $e_{k+1}^{\top}
    = e_k^{\top} M + e_1^{\top}$. The term $e_k^{\top}
    = e_1^{\top}\sum_{i=0}^{k-1} M^i$ fulfils the recurrent relation.

    \textit{ii)}: For every matrix $M\in Sl(\N,2)$ there exists
    $k\in\{2,3,4\}$ such that $\sum_{i=0}^{k-1}M^i$ has all elements
    even, i.e., $\sum_{i=0}^{k-1}M^i\bmod 2 =
    \left(\begin{smallmatrix}
        0 & 0 \\ 0 & 0
    \end{smallmatrix}\right)$. Indeed, every matrix $M\in Sl(\N,2)$ is equal $\bmod\ 2$
    to one of the following 6 matrices
    \begin{equation}\label{matrices_classes}
        \begin{aligned}
        M_1 &= \begin{pmatrix}
            1 & 0 \\ 0 & 1
        \end{pmatrix}, &
        M_2 &= \begin{pmatrix}
            0 & 1 \\ 1 & 0
        \end{pmatrix}, &
        M_3 &= \begin{pmatrix}
            1 & 0 \\ 1 & 1
        \end{pmatrix}, \\
        M_4 &= \begin{pmatrix}
            1 & 1 \\ 0 & 1
        \end{pmatrix}, &
        M_5 &= \begin{pmatrix}
            1 & 1 \\ 1 & 0
        \end{pmatrix}, &
        M_6 &= \begin{pmatrix}
            0 & 1 \\ 1 & 1
        \end{pmatrix}.
        \end{aligned}
    \end{equation}
    By inspection of these matrices we get that $k=2$ if $M=M_1 \bmod 2$,
    $k = 4$ if $M=M_i \bmod 2$ for $i\in\{2,3,4\}$ and $k = 3$
    if $M=M_i \bmod 2$ for $i\in\{5,6\}$. E.g.
    $$ M_6^0 + M_6^1 + M_6^2 =
    \begin{pmatrix}
        1 & 0 \\ 0 & 1
    \end{pmatrix} +
    \begin{pmatrix}
        0 & 1 \\ 1 & 1
    \end{pmatrix} +
    \begin{pmatrix}
        1 & 1 \\ 1 & 0
    \end{pmatrix}  =
    \begin{pmatrix}
        0 & 0 \\ 0 & 0
    \end{pmatrix} \bmod 2.$$

    \textit{iii)}: For $k\in\{2,3,4\}$ associated to the matrix
    $M=M_\varphi$ by \textit{ii)} it holds that
    $$P Q^k P^{-1} \in Sl(\Z,3).$$
    Indeed,

    \[
        PQ^k P^{-1}
        = P
        \begin{pmatrix}
            M^k & \vec{0} \\ e_k^{\top} & 1
        \end{pmatrix} P^{-1}
        =
        \begin{pmatrix}
            M^k & \vec{0} \\ f_k^{\top} & 1
        \end{pmatrix}
    \]
    where
    \begin{equation*}
        \begin{aligned}
            f_k^\top
            = (\tfrac{1}{2},0)M^k + \tfrac{1}{2}\left(e_k^\top - (1,0)\right)
            &= \tfrac{1}{2} \left((1,0)(M^k - I) + e_1^\top \sum_{i=0}^{k-1}M^i\right)\\
            &= \tfrac{1}{2} \left((1,0)(M - I) + e_1^\top\right) \sum_{i=0}^{k-1}M^i.
        \end{aligned}
    \end{equation*}
    As by \textit{ii)} the matrix $\sum_{i=0}^{k-1}M^i$ has all
    elements even, $f_k^\top \in \Z^2$
    and thus $PQ^k P^{-1}\in\Z^{3\times 3}$.
    Moreover, $\det PQ^k P^{-1} = \det M^k =1$.

    \textit{iv)}: We prove that $P Q^k P^{-1} \in \mathcal{E}$.
    We denote $Q^{k}= \left(\begin{smallmatrix} \widetilde{A} &  \widetilde{B}& 0 \\ \widetilde{C} &
    \widetilde{D} & 0 \\ \widetilde{E} & \widetilde{F} & 1 \end{smallmatrix}\right)$. As $Q^k = \Rcal(\varphi^k)\in\mathcal{E}$,
    by Lemma~\ref{inequalities} we have that $\widetilde{A},
    \widetilde{B},\widetilde{C},\widetilde{D},\widetilde{E},\widetilde{F}
    \in\N$ and the inequalities
    \begin{equation}\label{eq:inequalities}
         \widetilde A \widetilde D - \widetilde B \widetilde C = 1, \quad
         \widetilde{E}<\widetilde{A}+\widetilde{C}, \quad
         \widetilde{F}<\widetilde{B}+\widetilde{D}, \quad
         -\widetilde C \leq \widetilde C \widetilde F - \widetilde D \widetilde E < \widetilde D.
    \end{equation}
    Analogous inequalities for
    $$ PQ^k P^{-1}
    = P \left(\begin{smallmatrix} \widetilde{A} &  \widetilde{B}& 0 \\ \widetilde{C} &
    \widetilde{D} & 0 \\ \widetilde{E} & \widetilde{F} & 1 \end{smallmatrix}\right) P^{-1}
    = \left(\begin{smallmatrix} \widetilde{A} &  \widetilde{B}& 0 \\ \widetilde{C} &
    \widetilde{D} & 0 \\
    \tfrac{\widetilde{E} + \widetilde{A} - 1}{2}
    & \tfrac{\widetilde{B} + \widetilde{F}}{2}  & 1 \end{smallmatrix}\right),$$
    namely $\tfrac{\widetilde{E} + \widetilde{A} - 1}{2}<\widetilde{A}+\widetilde{C}$,
    $\tfrac{\widetilde{B} + \widetilde{F}}{2}<\widetilde{B}+\widetilde{D}$
    and $-\widetilde C \leq \widetilde C \tfrac{\widetilde{B} + \widetilde{F}}{2} - \widetilde D \tfrac{\widetilde{E}
    + \widetilde{A} - 1}{2} < \widetilde D$,
    are a consequence of \eqref{eq:inequalities}.
    Moreover, $\tfrac{\widetilde{E} + \widetilde{A} - 1}{2}\in\N$
    as $\widetilde{A},\widetilde{E}\in\N$ and
    $\tfrac{\widetilde{E} + \widetilde{A} - 1}{2}\in\Z$ by \textit{iii)}.
    The result follows by Lemma~\ref{inequalities}.
    \end{proof}

If $\uu$ is a characteristic  Sturmian sequence, we provide an  algorithm finding the morphism fixing $\sqrt{\uu}$.    To state our result  we need to recall that
a word $w = w_0w_1 \dots w_{n-1}$ is a \emph{palindrome} if it reads the same from the left as from the right, i.e., $w_k = w_{n-1-k}$ for each $k=0, 1, \ldots, n-1$.
\begin{corollary}\label{naseodmocniny} Let $\uu\in\{0,1\}^\mathbb{N}$ be a  characteristic Sturmian sequence
    fixed by a primitive morphism $\varphi\in\M$ having the incidence matrix  $M=M_\varphi$. Let $k$ be the smallest positive integer such that $(1,1)M^k = (1,1)\bmod 2$. Then  $k\leq 3$ and the square root $\sqrt{\uu}$ is fixed by
    a morphism $\psi\in\M$ which is a conjugate of $
    \varphi^k$.   Moreover, $\psi(0)$ and $\psi(1)$ are
palindromes of odd length.  \end{corollary}
\begin{proof}
We start by showing that in this case it suffices to consider $k\leq 3$ in the proof of \Cref{naseodmocniny2}.
The characteristic sequence has a vector of parameters $\vec{v}(\uu)=
    (1-\alpha, \alpha, \alpha)^\top$, which is an eigenvector of the matrix  $\mathcal{R}(\varphi) =\begin{pmatrix}
    A & B & 0 \\ C & D & 0 \\ E & F & 1 \end{pmatrix}$.  From the relation
    $\Rcal (\varphi)\vec{v}(\uu)=\Lambda \vec{v}(\uu)$ we obtain the equality
    $$C(1-\alpha)+D\alpha = \Lambda\alpha = E(1-\alpha)+F\alpha + \alpha.$$
    As the values $\alpha$ and $1-\alpha$ are linearly independent over $\Q$, we conclude that
    $E=C$ and $F=D-1$. Hence in the proof of the previous theorem $e_1^\top = (E,F)=(C,D-1) = (0,1)(M -I)$ and consequently, $f_k^\top = \frac12 (1,1)\bigl(M^k - I)$.
    The assumption $(1,1)M^k = (1,1)\bmod 2$ is equivalent to the claim that $f_k^\top$ has integer coordinates.  One can easily check by inspection of  6 matrices listed in \eqref{matrices_classes} that the vector $f_k^\top$ has integer coefficients for some $k \in \{1,2,3\}$.   In details,
     \begin{enumerate}
        \item If $M  = M_i \bmod 2$ for $i\in\{1,2\}$,
then    $(1,1)(M_i-I) \bmod 2 = (0,0)$ and we take $k=1$.
        \item If $M  = M_i\bmod 2$ for $i\in\{3,4\}$,  then $(1,1)(M^{2}_i-I) \bmod 2 = (0,0)$,   but  $(1,1)(M_i-I) \bmod 2 \neq  (0,0)$. We take $k=2$.
        \item If $M  = M_i\bmod 2$ for $i\in\{5,6\}$,  then  $k=3$ is the smallest positive $k \in \mathbb{N}$ satisfying  $(1,1)(M^{k}_i-I) = (0,0)\bmod 2 $.
    \end{enumerate}

    \bigskip
 Finally let us show that $\psi(0)$ and $\psi(1)$ are palindromes. Let us recall
    that $\sqrt{\uu}={\bf s}_{\alpha,\frac{1}{2}}$.   Consider the biinfinite sequence
    ${\bf s}_{\alpha,\frac{1}{2}} =  \dots
    \nu_{-3}\nu_{-2}\nu_{-1}\nu_{0}\nu_{1}\nu_{2}\dots$.  Let us deduce that  its left part
     $ \dots
    \nu_{-3}\nu_{-2}\nu_{-1}$ is the mirror image of the right part $ \nu_{0}\nu_{1}\nu_{2}\dots$.

    As ${\bf s}_{\alpha,\frac{1}{2}}(n) = \lfloor\alpha(n+1) +\frac{1}{2}\rfloor-\lfloor
    \alpha n +\frac{1}{2}\rfloor$ for each $n\in\Z$, we have
    ${\bf s}_{\alpha,\frac{1}{2}}(-n-1)=\lfloor\alpha(-n) +\frac{1}{2}\rfloor-\lfloor
    \alpha (-n-1) +\frac{1}{2}\rfloor$. Using the relations
    $\lfloor -x \rfloor = -\lfloor x \rfloor - 1$
    and $\lfloor x + 1 \rfloor = \lfloor x \rfloor + 1$ for  $x\notin\Z$
    we get that ${\bf s}_{\alpha,\frac{1}{2}}(-n-1) = {\bf s}_{\alpha,\frac{1}{2}}(n)$
    for each $n\in\Z$.  It confirms  the  mirror symmetry.   In particular we have that $\nu_{-1}=\nu_{0}$ and hence $\psi(\nu_0)=\psi(\nu_{-1})$.
    Now we use the result of \cite{BaMaPe}, where  the authors proved that  a biinfinite Sturmian   sequence ${\bf s}_{\alpha, \delta}$ is fixed by a primitive morphism  $\psi \in \mathcal{M}$ if and only if its right part is fixed by  $\psi$.

The symmetry of
\begin{align*}
 \dots \nu_{-3}\nu_{-2}\nu_{-1} &= \dots\psi(\nu_{-3})\psi(\nu_{-2})\psi(\nu_{-1}) \quad \text{ and } \\ \nu_0\nu_1\nu_2 \dots &= \psi(\nu_{0})\psi(\nu_{1})\psi(\nu_{2})\dots
\end{align*}
gives that the mirror image of $\psi(\nu_{-1})$ is equal to $\psi(\nu_0)$.
    In other words $\psi(\nu_0)$ is a palindrome. The mirror symmetry also
    gives that the image under $\psi$ of the letter other than $\nu_0$
    is a palindrome as well.

 The morphism $\psi$ has the incidence matrix $M^k$. The lengths of images of the letters $0$ and $1$  by $\psi$ fulfil $(|\psi(0)|, |\psi(1)|) = (1,1)M^k$. From the assumption $(1,1)M^k =(1,1)\bmod 2$ we have that the palindromes $\psi(0)$ and $\psi(1)$ have odd length.
\end{proof}

\begin{example} (continuation of Example \ref{priklad})\ \
    We observe that $\varphi = DG^2$
    which fixes a characteristic Sturmian sequence has the incidence matrix
    $M = \left(\begin{smallmatrix} 1 & 2 \\ 1 & 3 \end{smallmatrix}\right)$ and  $(1,1)\bigl(M^2 - I\bigr) = (0,0)\bmod 2$.
    The morphism $\varphi^{2}$ has  a conjugate
    $$\psi:0\mapsto 1010101 , 1\mapsto 1010101101011010101, $$
    such that  both  $\psi(0)$ and $\psi(1)$ are palindromes of odd length.
    This corresponds with Corollary~\ref{naseodmocniny}.
    The  prefix of $\sqrt{\uu}$ displayed in Example \ref{priklad} illustrates that $\sqrt{\uu}$ is fixed by $\psi$.

    To illustrate  Corollary~\ref{naseodmocniny}
    more thoroughly, we show other examples
    where a conjugate of $\varphi$, $\varphi^2$ or $\varphi^3$
    is used to fix the square root $\sqrt{\uu}$.
    \begin{table}[h]
    \centering
    \begin{tabular}{lcl}
    $\varphi~~~~~~~~~~$ & $k$ & \qquad conjugate of $\varphi^k$  \\\hline
    $D^2G^2$  & $1$ & \qquad $0\mapsto 101$, $1\mapsto 1011101$  \\
    $GDG$  & $1$ & \qquad $0\mapsto 010$, $1\mapsto 01010$  \\
    $D^2G$  & $2$ & \qquad $0\mapsto 10111011101$, $1\mapsto 101110111011101$ \\
    $DG$  & $3$ & \qquad $0\mapsto1010110110101, 1\mapsto101011011010110110101$\\
    $GD$  & $3$ & \qquad $0\mapsto010100100101001001010, 1\mapsto0101001001010$
\end{tabular}
\end{table}

\end{example}

We show an example of a morphism $\varphi$
for which only a conjugate of $\varphi^4$
fixes $\sqrt{\uu}$. Of course, the fixed point $\uu$ of $\varphi$ is not  characteristic.
\begin{example}
    The morphism $\varphi=DG\widetilde{G}$ has the representation
    $\Rcal(\varphi) = R_D R_G R_{\widetilde{G}} = \left(\begin{smallmatrix} 1 & 2 & 0\\ 1 & 3 & 0 \\ 1 & 3 & 1\end{smallmatrix}\right)$,
    the incidence matrix
    $M = \left(\begin{smallmatrix} 1 & 2 \\ 1 & 3 \end{smallmatrix}\right)$
    and $e_1^\top = (1,3)$.
    In the proof of Theorem~\ref{naseodmocniny2}
    $f_k^\top = \tfrac{1}{2} \left((1,0)(M - I) + e_1^\top\right) \sum_{i=0}^{k-1}M^i
    \in\Z^2$ if either $\sum_{i=0}^{k-1}M^i$ has all elements even
    or $\left((1,0)(M - I) + e_1^\top\right)$ has all elements even.
    As $\left((1,0)(M - I) + e_1^\top\right) = (1,5)$,
    the smallest possible $k\in\{1,2,3,4\}$ such that a conjugate
    of $\varphi^k$ fixes $\sqrt{\uu}$ is $k=4$.
\end{example}


\section{Comments}
Let us formulate two questions on Sturmian morphisms we were not able to answer.

\begin{itemize}
\item Conjugate Sturmian morphisms can be ordered with respect to the relation $\triangleright$ from the leftmost conjugate morphism to the rightmost conjugate.  All the morphisms have the same incidence matrix  and thus the matrices of their faithful representation differ only in the last row. It will be interesting to describe the relationship between the last rows of two subsequent elements in this chain.   Let us demonstrate the question on an example.

By Theorem \ref{sturmian_conjugate}, $M=\left( \begin{array}{ll}
1&2\\
2&5
\end{array}\right) $
 is the incidence matrix for 9 mutually conjugate Sturmian morphisms. Let us list them  from the leftmost conjugate to the rightmost conjugate morphism:
 $$ \psi_0 =  \widetilde{D}^2\widetilde{G}^2,  \ \ \psi_1 = \widetilde{D}D\widetilde{G}^2, \ \ \ \psi_2={D}^2G\widetilde{G}^2, \ \ \ \psi_3 = \widetilde{D}^2G\widetilde{G}, \ \ \ \psi_4 = D\widetilde{D}G\widetilde{G}$$
 $$\psi_5 = {D}^2G\widetilde{G}, \ \ \ \psi_6 = \widetilde{D}^2{G}^2, \ \ \ \psi_7 = \widetilde{D} D{G}^2, \  \ and \ \ \psi_8 = {D}^2{G}^2
. $$
If we compute the faithful representation of the morphism $\psi_i$ for $i=0,1,\ldots, 8$, we obtain   $\mathcal{R}(\psi_i)= \left( \begin{array}{ll}
M&0\\
e_i&1
\end{array}\right)$ with
$
e_0 = (0,2), \ e_1 = (1,4), \ e_2=(2,6), \ e_3 = (0,1)$, $ e_4 = (1,3), \ e_5 = (2,5), \ e_6 = (0,0), \ e_7 = (1,2), \ e_8 = (2,4).
$

The question is how to generate the sequence $e_i$ using the form of the matrix $M$, and, consequently, the sequence $\rho_i$, where $\rho_i$ is the intercept of the fixed point of~$\psi_i$.

\item \Cref{Yasutomi_necessary} is only the necessary condition of Yasutomi's theorem. The  sufficient condition of the theorem can be formulated in our formalism as follows:

Let $x,y,z \in \mathbb{Q}(\sqrt{m})$, where $m \in \mathbb{N}$,   $\sqrt{m}\notin \mathbb{Q}$ and $\frac{x}{y} \notin \mathbb{Q}$.  If   $ (x,y,z)^\top \in C_1$ and the Galois conjugate  $(\overline{x},\overline{y},\overline{z})^\top \in C_2$, then  $ (x,y,z)^\top $  is an eigenvector of a matrix  $R \in \mathcal{E}$, $R\neq I_3$.

 We were not able to prove this statement directly and obtain a proof that would differ substantially from Yasutomi's proof.
\end{itemize}

The exhibited faithful representation relies on the geometric representation of  Sturmian sequences. A  natural question is finding a faithful representation of some other monoids of morphisms that fix some other families of infinite sequences, e.g., Arnoux--Rauzy sequences and sequences coding symmetric $k$-interval exchange transformations.
Let us point out that a geometrical representation of  ternary Arnoux-Rauzy sequences can be found  already in the article \cite{ArRo}.
Recently, a representation of Arnoux--Rauzy sequences using a generalized  Ostrowski numeration system is introduced in \cite{Pelto21}.
The sequences coding $3$-interval exchanges can be viewed as sequences obtained by cut-and-project sets, see \cite{MaPeSta1}.

\section*{Acknowledgements}
Jana Lep\v{s}ov\'{a} acknowledges financial support by The French Institute in Prague and the Czech Ministry of Education, Youth and Sports through the Barrande fellowship programme and Agence Nationale de la Recherche through the project Codys (ANR-18-CE40-0007), and the support by  Grant Agency of Czech technical university in Prague, through the project   SGS20/183/OHK4/3T/14.
Edita Pelantov\'{a} acknowledges financial support by The Ministry of Education, Youth and Sports of the Czech Republic, project no. {CZ.02.1.01/0.0/0.0/16\_019/0000778}.
\v{S}t\v{e}p\'{a}n Starosta acknowledges the support of the OP VVV MEYS funded
project {CZ.02.1.01/0.0/0.0/16\_019/0000765}.

\bibliographystyle{siam}
\IfFileExists{biblio.bib}{\bibliography{biblio}}{\bibliography{../!bibliography/biblio}}

\begin{thebibliography}{10}

\bibitem{ArRo}
{\sc P.~Arnoux and G.~Rauzy}, {\em Repr\'esentation g\'eom\'etrique de suites
  de complexit\'e $2n+1$}, Bull. Soc. Math. France, 119 (1991), pp.~199--215.

\bibitem{BaMaPe}
{\sc P.~Bal\'{a}\v{z}i, Z.~Mas\'{a}kov\'{a}, and E.~Pelantov\'{a}}, {\em
  Complete characterization of substitution invariant {S}turmian sequences},
  Integers, 5 (2005), pp.~A14, 23.

\bibitem{BEIR}
{\sc V.~Berth\'{e}, H.~Ei, S.~Ito, and H.~Rao}, {\em On substitution invariant
  {S}turmian words: an application of {R}auzy fractals}, Theor. Inform. Appl.,
  41 (2007), pp.~329--349.

\bibitem{Dekking}
{\sc M.~Dekking}, {\em Substitution invariant {S}turmian words and binary
  trees}, Integers, 18A (2018), p.~\#A17.

\bibitem{DuSch}
{\sc N.~Dunford and J.~T. Schwartz}, {\em Linear operators. {P}art {I}}, Wiley
  Classics Library, John Wiley \& Sons, Inc., New York, 1988.
\newblock General theory, With the assistance of William G. Bade and Robert G.
  Bartle, Reprint of the 1958 original, A Wiley-Interscience Publication.

\bibitem{Fiedler}
{\sc M.~Fiedler}, {\em Special matrices and their applications in numerical
  mathematics}, Dover Publications, Inc., Mineola, NY, second~ed., 2008.
\newblock Translated from the Czech by Petr P\v{r}ikryl and Karel Segeth.

\bibitem{Fogg}
{\sc N.~P. Fogg}, {\em Substitutions in Arithmetics, Dynamics and
  Combinatorics}, vol.~1794 of Lecture notes in mathematics, Springer, 1st~ed.,
  2002.

\bibitem{algebra}
{\sc H.~Georgi}, {\em Lie Algebras In Particle Physics: from Isospin To Unified
  Theories}, CRC Press, 2018.

\bibitem{Lo2}
{\sc M.~Lothaire}, {\em Algebraic Combinatorics on Words}, no.~90 in
  Encyclopedia of Mathematics and its Applications, Cambridge University Press,
  2002.

\bibitem{MaPeSta1}
{\sc Z.~Masáková, E.~Pelantová, and {\v{S}}.~Starosta}, {\em Exchange of
  three intervals: substitutions and palindromicity}, Eur. J. Combin., 62
  (2017), pp.~217--231.

\bibitem{MiPa_Rauzy}
{\sc F.~Mignosi and P.~S\'e\'ebold}, {\em Morphismes sturmiens et r\`egles de
  {Rauzy}}, J. Th\'eor. Nombres Bordeaux, 5 (1993), pp.~221--233.

\bibitem{HMo}
{\sc M.~Morse and G.~A. Hedlund}, {\em Symbolic dynamics {II}. {S}turmian
  trajectories}, Amer. J. Math., 62 (1940), pp.~1--42.

\bibitem{PeWh}
{\sc J.~Peltom\"{a}ki and M.~A. Whiteland}, {\em A square root map on
  {S}turmian words}, Electron. J. Combin., 24 (2017), pp.~Paper No. 1.54, 50.

\bibitem{Pelto21}
{\sc J.~Peltomäki}, {\em Initial nonrepetitive complexity of regular
  episturmian words and their {D}iophantine exponents}.
\newblock preprint available at \url{
  https://doi.org/10.48550/arXiv.2103.08351}, 2021.

\bibitem{Peng}
{\sc L.~Peng and B.~Tan}, {\em Sturmian sequences and invertible
  substitutions}, Discr. Math. Theoret. Comput. Sci., 13 (2011), pp.~63--68.

\bibitem{Sa}
{\sc K.~Saari}, {\em Everywhere {$\alpha$}-repetitive sequences and {S}turmian
  words}, European J. Combin., 31 (2010), pp.~177--192.

\bibitem{Ya99}
{\sc S.-I. {Yasutomi}}, {\em {On {S}turmian sequences which are invariant under
  some substitutions.}}, in {Number theory and its applications. Proceedings of
  the conference held at the RIMS, Kyoto, Japan, November 10--14, 1997}, vol.~2
  of Dev. Math., Dordrecht: Kluwer Academic Publishers, 1999, pp.~347--373.

\end{thebibliography}

\end{document}